\documentclass[twoside,11pt,a4paper]{article}
\usepackage{amsmath,float,amssymb,mathtools,amsthm,mathrsfs,multirow,booktabs,setspace,morefloats,epsfig,caption,graphicx,inputenc,subcaption}
\usepackage{bigints}
\usepackage{xcolor}
\usepackage{placeins}
\usepackage{titling}
\newtheorem{lemma}{Lemma}
\newtheorem{thm}{Theorem}
\newtheorem{cor}{Corollary}
\usepackage{fancyhdr}
\newcommand\shorttitle{Some parametric tests based on sample spacings}

\newcommand\authors{Rahul Singh and Neeraj Misra}

\allowdisplaybreaks
\title{Some parametric tests based on sample spacings}

\author{Rahul Singh$^1$ and Neeraj Misra$^2$}
\date{
	\small{Department of Mathematics and Statistics, Indian Institute of Technology Kanpur, India}\\%
}

\begin{document}
	\maketitle
	{\let\thefootnote\relax\footnote{{ Email: $^1$sirahul@iitk.ac.in,  $^2$neeraj@iitk.ac.in}}}
\begin{abstract}
	Assume that we have a random sample from an absolutely continuous distribution (univariate, or multivariate) with a known functional form and some unknown parameters. In this paper, we have studied several parametric tests based on statistics that are symmetric functions of $m$-step disjoint sample spacings. Asymptotic properties of these tests have been investigated under the simple null hypothesis and under a sequence of local alternatives converging to the null hypothesis. The asymptotic properties of the proposed tests have also been studied under the composite null hypothesis. We observed that these tests have similar asymptotic properties as the likelihood ratio test. Finite sample performances of the proposed tests are assessed numerically. A data analysis based on real data is also reported. The proposed tests provide alternative to similar tests based on simple spacings (i.e., $m=1$), that were proposed earlier in the literature. These tests also provide an alternative to  likelihood ratio tests in situations where likelihood function may be unbounded and hence, likelihood ratio tests do not exist.\\~\\
 \textit{Keywords: }{Asymptotic distribution, generalised spacings estimator, hypothesis test,  likelihood ratio test, multivariate spacings, nearest neighbour, sample spacings.}
\end{abstract}
\section{Introduction}
\label{intro}
Let $X_1,X_2,\ldots,X_n$ be independently and identically distributed (iid) random variables having an absolutely continuous distribution function $F_\eta$ with $\eta\in \Theta\subseteq \mathbb{R}^p$. Assume that for every $\theta\in\Theta$, the functional form of $F_\theta$ is specified and that the true parameter $\eta\in\Theta$ is unknown. Here, we are interested in tests for simple and composite hypotheses concerning the unknown true parameter. There are many ways to address this problem. A popular method is to develop a test statistic based on sample spacings. Let $X_{1:n},X_{2:n},\ldots,X_{n:n}$ denote the order statistics corresponding to $X_1,X_2,\ldots,X_{n}$. Define $X_{0:n}=-\infty$ and $X_{n+1:n}=\infty$.  For any positive integer $m\ (<n)$, the $m$-step disjoint sample spacings are defined as
\begin{align*}
	D^{(m)}_{j,n}(\theta)= F_\theta(X_{jm:n})-F_\theta(X_{(j-1)m:n})\text{ for }\ j=1,2,\ldots, \left\lfloor{\frac{n+1}{m}} \right\rfloor, 
\end{align*}
where for any real number $x$, $\lfloor{x} \rfloor$ denotes the largest integer not exceeding $x$.
For any positive integer $m~(=m(n))$, sufficiently smaller than $n$, $ \left\lfloor{\dfrac{n+1}{m}} \right\rfloor\simeq\ \dfrac{n+1}{m}=M$, (say). For stating asymptotic results, without loss of generality, we can take $M$ to be an integer. 
Let $\theta_0\in\Theta$ be pre-specified. For testing $H_0:\eta=\theta_0$, a useful test statistic based on disjoint sample spacings is of the form 
\begin{align}
	S_{\phi,n}^{(m)}(\theta)=\frac{1}{M}\sum_{j=1}^{M}\phi\pmb(MD_{j,n}^{(m)}(\theta)\pmb),\ \theta\in\Theta, \label{eq2}
\end{align}
for some real-valued (convex, or concave) function $\phi$ defined on the positive half of the real line. The choice of function $\phi$ corresponding to some popular test statistics are as follows:
\begin{table*}[h!]
	\begin{center}
		\begin{tabular}{l l} 
			\hline
			$\phi(x)$ & Statistic \\
			\hline 
			$x^2$ & Greenwood Statistic (Greenwood 1946) \\
			
			$\log(x)$ & Log Spacing Statistic (Moran 1951)\\
			
			$|x-1|^r, r>0$ & Generalized Rao's Spacing Statistic (Rao 1969)\\
			
			$x^r,\ r> 0$ & Kimball Statistic (Kimball 1974)\\
			
			$x\log( x)$ & Relative Entropy Spacing Statistic (Misra and van der Meulen 2001)\\
			\hline
		\end{tabular}
	\end{center}
\end{table*}
\FloatBarrier
Ekstr\"om (2013) studied parametric tests based on statistics of the type (\ref{eq2}) with simple spacings (i.e., $m=1$). The goal of this article is to extend the results of Ekstr\"om (2013) to test statistics of type (\ref{eq2}), based on $m$-step disjoint spacings. We also extend the study of Ekstr\"om (2013)  to the multivariate setup. In the sequel, we describe some advantages of considering test statistics based on higher order disjoint spacings. Data sets arising from various real life situations may contain ties (due to rounding off, etc.). In such situations, the optimal test (one corresponding to $\phi(x)=-\log(x)$) proposed by Ekstr\"om (2013) can not be used and a remedy may be obtained by using tests based on higher order spacings. Further, in the context of testing goodness of fit (Del Pino 1979) and estimation (Ekstr\"om et al. 2020) the statistics based on higher order spacings are known to be more efficient than those based on simple spacings. The above discussion poses a natural question that whether the use of higher order disjoint spacings will also be beneficial for parametric testing problems. Thus, study of parametric tests based on higher order disjoint spacings is of theoretical and practical interest. It is worth mentioning here that Ekstr\"om et al. (2020) studied properties of estimators based on higher order disjoint spacings. Now, we will discuss some of the relevant studies on treatments of testing and estimation problems based on sample spacings and the relation of our work to these studies.

Goodness of fit tests, based on class of statistics defined by (\ref{eq2}), have been widely studied in the literature. For $m=1$ (simple spacings case), under quite general conditions on the underlying distributions and the function $\phi$, Sethuraman and Rao (1970) established the asymptotic normality of the test statistic $S_{\phi,n}^{(1)}(\theta_0)$ under the simple null hypothesis $H_0:\eta=\theta_0$. Del Pino (1979) extended this result by establishing the asymptotic normality of $S_{\phi,n}^{(m)}(\theta_0)$, under $H_0$, for any finite $m$. Mirakhmedov (2005) further extended these results to situations where $m$ is allowed to grow with $n$ such that $m=o(n)$. 
Goodness of fit tests based on test statistics in (\ref{eq2}) can detect only those alternatives that converge to the null distribution at the rate of $n^{-1/4}$ or slower. Within such class of statistics, the Greenwood test statistic is known to be asymptotically locally most powerful in terms of Pitman efficiency (see Sethuraman and Rao 1970). 

Cheng and Amin (1983), and Ranneby (1984) studied estimation of the true unknown parameter $\eta$ and proposed maximising the parametric function (\ref{eq2}) with $\phi(x)=\log(x),\ x>0$. The estimator so obtained is called the maximum spacings product estimator (MSPE). Ghosh and Jammalmadaka (2001) continued this study further by considering a general $\phi$. Under quite general conditions, they showed that such an estimator has similar asymptotic properties as the maximum likelihood estimator (MLE). Ekstr\"om et al. (2020) extended the results of Ghosh and Jammalmadaka (2001) to situations where $m$ is any finite, but fixed, positive integer or $m\to\infty$ such that $m=o(n)$. Such an estimator is known as the generalised spacing estimator (GSE).

Note that $S_{\phi,n}^{(m)}(\theta_0)$ alone can not be used as a test statistic as it does not contain information regarding alternative hypothesis. Suppose that $\hat{\eta}$ is a $\sqrt{n}$-consistent estimator of $\eta$, then, for large $n$,  $S_{\phi,n}^{(m)}(\hat{\eta})$ is expected to be close to $S_{\phi,n}^{(m)}(\theta_0)$, under $H_0:\eta=\theta_{0}$. Hence, some distance function that measures departure of $S_{\phi,n}^{(m)}(\theta_0)$ from $S_{\phi,n}^{(m)}(\hat{\eta})$ can be used as a test statistic, where a large value of the distance would indicate incompatibility of the data with $H_0$. Based on this idea, 
Torabi (2006) proposed parametric tests based on simple spacings. Ekstr\"om (2013) further studied such parametric tests based on simple spacings and showed that such tests have asymptotic properties similar to the likelihood ratio test (LRT). A part of this paper can be seen as an extension of work of Ekstr\"om (2013) to statistics based on general $m$-step disjoint spacings.

In the multivariate setup, Zhou and Jammalamadaka (1993) generalised the concept of univariate spacings using nearest neighbour balls. Kuljus and Ranneby (2015) proved consistency of the GSE for multivariate observations. Under some conditions, Kuljus and Ranneby (2020) found that the GSE for multivariate observations have similar asymptotic properties as the MLE. In this paper we have studied parametric tests based on multivariate sample spacings and found that the corresponding test statistics have similar asymptotic properties as their univariate analogues.  {In other words, asymptotically, the test statistics based on multivariate sample spacings are  distributed as central chi-square under $H_0$, and non-central chi-square under a sequence of local alternatives.}

We have also performed an extensive numerical study to assess finite sample performances of the proposed tests. For large sample sizes ($n\geq 50$), asymptotic distributions of the test statistics are found to be reasonably close to their limiting distribution for many models that we have investigated. It is noteworthy that the proposed tests can be used in situations where LRT may not exist due to unboundedness of the likelihood function, e.g., certain finite mixture distributions and heavy tailed distributions (see Pitman 2018, p.~70). Even in certain situations where LRT exist, performances of proposed tests are observed to be comparable to those of LRT.

Throughout, $\chi^2_d(\delta)$ denotes a non-central chi square random variable with $d$ degrees of freedom and non-centrality parameter $\delta$, and $\chi_d^2$ denotes a central chi square random variable with $d$ degrees of freedom. $\mathbb{R}$ will denote the real line, and for any positive integer $p$, $\mathbb{R}^p$ will denote the $p$-dimensional Euclidean space.  Convergences in probability and distribution are denoted by $\xrightarrow{p}$ and $\xrightarrow{d}$, respectively.

The rest of the paper is structured as follows. Tests for univariate distributions are introduced in Section 2. In Section 3, tests for multivariate distributions are discussed. A numerical study to assess finite sample performance of the proposed tests is reported in Section 4. In Section 5, a real data analysis is presented. A discussion based on the study is made in Section 6. All the proofs are given in the Appendix.
\section{Tests for Univariate Distributions}
Let $\eta\in\Theta$ be the true and unknown value of the parameter. Based on a random sample $X_1,X_2,\ldots,X_n$ from $F_\eta$, we aim to test the following hypothesis
\begin{align}\label{simpleh0}
	H_0:\eta=\theta_0 \text{ against } H_A:\eta\neq\theta_0, 
\end{align}
where $\theta_0\in\Theta$ is pre-specified. Let $\phi:(0,\infty)\to\mathbb{R}$ be a convex function. As a measure of departure of any $F_\theta,\ \theta\in\Theta$ from $F_\eta$, Csisz\'ar (1977) defined the $\phi$-divergence of $F_{\theta}$ with respect to $F_{\eta}$ as
\begin{align*}
	S_\phi (\theta,\eta):=\int_{-\infty}^{\infty}\phi\left(\frac{f_{\theta}(x)}{f_{\eta}(x)}\right)f_{\eta}(x)\ dx,\ \theta\in\Theta,
\end{align*}
where $f_\tau$ is the density function corresponding to the distribution function $F_\tau,\ \tau\in\Theta$. If $\phi(x)=-\log(x),\ x>0$, the $\phi$-divergence is known as the Kullback-Leibler (KL) divergence. Using Jensen's inequality, we have
\begin{align*}
	\inf_{\theta\in\Theta}S_\phi(\theta,\eta) = S_\phi(\eta,\eta)=\phi(1).
\end{align*}
Thus, if $\hat{\eta}$ is a suitable estimator of $\eta$, then 
\begin{align*}
	{T}_{\phi,n}(\theta_0)= S_\phi(\theta_0,\hat{\eta})-\inf_{\theta\in\Theta}S_\phi(\theta,\hat{\eta})
\end{align*}
can be used as a measure of departure from the null hypothesis ($H_0:\eta=\theta_0$). 
Since $\inf_{\theta\in\Theta}S_\phi(\theta,\eta) = S_\phi(\eta,\eta)$, a suitable estimator of $\eta$ is the one that minimises $S_\phi(\theta,\eta)$ with respect to $\theta\in \Theta$. As $S_\phi(\theta,\eta)$ involves unknown $\eta$, an appropriate approximation of $S_\phi(\theta,\eta)$ may be used for this purpose.  Note that $\dfrac{n}{m}\{F_\theta({X_{jm:n}})-F_\theta({X_{(j-1)m:n}})\}$ is a non-parametric histogram estimator of $\dfrac{f_\theta(x)}{f_{\eta}(x)}$, $x\in\ {[X_{(j-1)m:n},\ X_{jm:n})}$ (see Prakasa Rao 1983, Section 2.4, pp.~102-115). 
This suggests that $S_{\phi,n}^{(m)}(\theta)$, defined in (\ref{eq2}), is an appropriate approximation of $S_\phi(\theta,\eta)$ and $\hat{\theta}_{\phi,n}^{(m)}=\arg\inf_{\theta\in\Theta}S_{\phi,n}^{(m)}(\theta)$ is a reasonable estimator of $\eta$. Thus, 
\begin{align*}
	T_{\phi,n}^{(m)}(\theta_0)= S_{\phi,n}^{(m)}(\theta_0)-\inf_{\theta\in\Theta}S_{\phi,n}^{(m)}(\theta)=S_{\phi,n}^{(m)}(\theta_0) -S_{\phi,n}^{(m)}(\hat{\theta}_{\phi,n}^{(m)})
\end{align*}
is a suitable statistic for testing $H_0:\eta=\theta_0$, where significantly large values of $T_{\phi,n}^{(m)}(\theta_0)$ provide an evidence of departure from the null hypothesis ($H_0:\eta=\theta_0$).

Ekstr\"om et al. (2020) assumed $m=o(n)$ and under quite general conditions, proved consistency and asymptotic normality of  $\hat{\theta}_{\phi,n}^{(m)}$ as an estimator of the true parameter $\eta$. They observed that the convergence is uniform, i.e.,
\begin{align}\label{eks2020}
	\lim_{n\to\infty}\sup_{x}\bigg|\mathbb{P}_\eta\left(\sqrt{n}(\hat{\theta}_{\phi,n}^{(m)}-\eta)\leq x\right)-\Phi\left(x\sqrt{\frac{I(\eta)}{\sigma_{\phi,m}^2}}\right) \bigg| 
	=0,
\end{align}
where
\begin{align*}
	&\sigma_{\phi,m}^2=\frac{mVar\left({\zeta}_{m} \phi '({\zeta}_{m} )\right)+(2m+1)\mu _{\phi,m}^{2} -2m\mu _{\phi,m} E\left({\zeta}_{m} ^{2} \phi '({\zeta}_{m} )\right)}{\left(E\left({\zeta}_{m}^{2} \phi ''({\zeta}_{m} )\right)\right)^2},\\ &{\zeta_m}\sim \frac{1}{m}Gamma(m,1),\ \mu_{\phi,m}= {E} \left({\zeta}_{m} \phi '({\zeta}_{m} )\right). 
\end{align*}
Here, $Gamma(m,1)$  denotes the gamma distribution with shape parameter $m$ and scale parameter $1$.

Earlier, the special case of $m=1$ was studied by Ghosh and Jammalmadaka (2001) and similar results were obtained. 
For $m=1$ and $\phi(x)=-\log(x)$, $\hat{\theta}_{\phi,n}^{(m)}$ was among the first studied estimators based on  sample spacings (see Cheng and Amin 1983; Ranneby 1984). Using the generalised mean value theorem, for some $\tilde{X}_j\in~{(X_{(j-1)m:n},\ X_{jm:n})}$, we have
\begin{align*}
	T_{-\log,n}^{(m)}(\theta_0)
	=& \sup_{\theta\in\Theta} \frac{1}{M}\sum_{j=1}^{M}\log\left(\frac{D_{j,n}^{(m)}({\theta})} {D_{j,n}^{(m)}({\theta_0})}\right)
	= \sup_{\theta\in\Theta} \frac{1}{M}\sum_{j=1}^{M}\log\left(\frac{f_{\theta}(\tilde{X}_j)} {f_{\theta_0}(\tilde{X}_j)}\right)\\
	=& \frac{1}{M}\log\left(  \frac{\sup_{\theta\in\Theta}\prod_{j=1}^{M}f_{\theta}(\tilde{X}_j)} {\prod_{j=1}^{M}f_{\theta_0}(\tilde{X}_j)}\right).
\end{align*}
Thus, the test based on $T_{-\log,n}^{(m)}(\theta_0)$ seems to be asymptotically equivalent to the LRT. 
The LRT is quite popular and are known to have nice asymptotic properties. Ekstr\"om (2013) showed that tests based on $T_{-\log,n}^{(1)}(\theta_0)$ have similar properties as the LRT. 

For the problem (\ref{simpleh0}), we propose to use the $m$-spacings based test statistic 
\begin{align}\label{teststt}
	\tilde{T}_{\phi,n}^{(m)}(\theta_0)=\frac{2n{T}_{\phi,n}^{(m)}(\theta_0)} {E\left({\zeta}_{m}^{2} \phi ''({\zeta}_{m} )\right)\sigma_{\phi,m}^2}. 
\end{align}
\subsection{Main Results}
For proving theoretical results of this paper, we will require the following assumptions:
\begin{itemize}
	\item[(A1) ] Distributions in the family $\{F_\theta,\ \theta\in\Theta\}$ have a common support and the true parameter $\eta$ is an interior point of $\Theta$;
	\item[(A2) ] Distributions in the family $\{F_\theta,\ \theta\in\Theta\}$ are identifiable i.e., if $\theta_1\neq\theta_2$, then $F_{\theta_1}(x)\neq F_{\theta_2}(x)$ for some $x\in \mathbb{R}$ and $F_\theta$ is differentiable with respect to $\theta\in\Theta$;
	\item[(A3) ] $\phi:(0,\infty)\to\mathbb{R}$ is a strictly convex and thrice continuously differentiable function. Moreover, $\text{Var}(Z_1\phi'(\zeta_m))$, $\mathbb{E}(Z_1^2\phi''(\zeta_m))$ and $\mathbb{E}(Z_1^3\phi'''(\zeta_m))$ are finite and bounded away from zero, where $Z_1,\ldots,Z_m$ are iid standard exponential variates and $\zeta_{m}=\frac{1}{m}(Z_1+\cdots+Z_m)$;
	\item[(A4) ] $f_\eta(x)$, $F_{\eta}^{-1} (x)$, $\dfrac{\partial}{\partial x}f_\eta(x)$ and $\dfrac{\partial^2 }{\partial x\ \partial \theta_j}f_\theta(x)\bigg|_{\theta=\eta} $ are continuous functions of $x$ for $j\in\{1,2,\ldots,p\}$. Moreover, for any $i,j,k\in\{1,2,\ldots,p\}$, $\dfrac{\partial}{\partial \theta_j}f_\theta(x)$, $\dfrac{\partial^2}{\partial\theta_j \partial\theta_k}f_\theta(x)$ and $\dfrac{\partial^3}{\partial\theta_j \partial\theta_k\partial\theta_l}f_\theta(x)$ are continuous functions of $x$ as well as $\theta$ in an open neighborhood of $\eta$ (say $\Theta_0$);
	\item[(A5) ] For any $j,k\in\{1,2,\ldots,p\}$,
	\begin{align*}
		\bigint_{0}^{1}\left(\frac{\dfrac{\partial^2}{\partial\theta_j \partial\theta_k}f_\theta(F_{\eta}^{-1}(u))\bigg|_{\theta=\eta}}{f_{\eta}(F_{\eta}^{-1}(u))}\right)\ du&< \infty,\\
		I_{j,k}(\theta)= \int_{-\infty}^{\infty}
		\frac{\frac{d}{d\theta_j}f_\theta(x)\frac{d}{d\theta_k}f_\theta(x)}{f_\theta(x)}\ dx &<\infty,\text{ for every } \theta \in\Theta_0,
	\end{align*} 
	and the Fisher information matrix $I(\theta)=(I_{j,k}(\theta))$ is positive definite, for every $\theta\in\Theta_0$;
	\item[(A6) ] $\lim_{t\to\infty}\dfrac{\min\{0,\phi(t)\}}{t}=0$ and $\big|\phi(t)\big|\leq a(t^{-b}+t^c)\ \forall t>0$, for some non-negative constants $a$, $b$ and $c$.
	\item[(A7) ]  {For $j,k,l=1,2,\ldots,p$, there exists a function $M_{jkl}$ such that
		\begin{align*}
			&\bigg|\frac{\partial^3}{\partial\theta_j\partial\theta_k\partial\theta_l} \phi\left(\frac{n+1}{m}D_{i:n}^{(m)}(\theta)\right)\bigg|\leq M_{jkl}\left(\frac{n+1}{m}D_{i:n}^{(m)}(\theta_0)\right),~\forall~\theta\in\Theta_0\\ &\text{ and }m_{jkl}=\mathbb{E}(M_{jkl}(\bar{\zeta}_m))<\infty.	
	\end{align*}} 
\end{itemize}
The following convex functions are of special interest, because these functions satisfy the conditions required for ensuring consistency and asymptotic normality of ${\theta}_{\phi,n}^{(m)}$, and $\lim_{m\to\infty}\sigma_{\phi,m}^2=1$ (cf. Ekstr\"om et al. 2020):
\begin{eqnarray}\label{pow_fun}
	\phi(x)=\phi_\gamma(x)=\left\{\begin{array}{ll}
		\gamma^{-1}(1+\gamma)^{-1}(x^{\gamma+1}-1), & \mbox{ if }\gamma\neq-1,0. \\
		-\log x, & \mbox{ if }\gamma=-1.\\
		x\log x, & \mbox{ if }\gamma=0.
	\end{array}\right. \label{eq5}
\end{eqnarray}

Under assumptions (A1), (A2) and (A6), Ekstr\"om et al. (2020) proved that  ${\theta}_{\phi,n}^{(m)}$ is a consistent estimator of the true parameter $\eta$. Ekstr\"om et al. (2020) commented that multi-parameter version of result (\ref{eks2020}) is possible but  they did not state the result explicitly. The multivariate version is desirable for our purpose, and we state it below.
\begin{thm}\label{thm1}
	Let $X_1,X_2,\ldots,X_{n}$ be iid observations from $F_{\eta}$.  Suppose that $m=o(n)$, and assume that conditions  {(A1)-(A7)} hold. Then,  {there exists a sequence $\{\hat{\theta}^{(m)}_{\phi,n}=$ 
		$\arg \inf_{\theta\in\Theta}S^{(m)}_{\phi,n}(\theta)\}_{n\geq1}$ such that $\hat{\theta}^{(m)}_{\phi,n}$ is a consistent estimator of $\eta$, and} 
	\begin{itemize}
		\item[(a)] for finite $m$,
		\begin{align*}
			\sqrt{n}({\hat{\theta}}_{\phi,n}^{(m)}-\eta) \stackrel{d}{\rightarrow} 
			N\left(0,\ {\sigma_{\phi,m}^2} {I(\eta)}^{-1}\right); 
		\end{align*}
		\item[(b)] for $m\to\infty$ and $\lim_{m\to\infty}\sigma_{\phi,m}^2=1$,  
		\begin{align*}
			\sqrt{n}({\hat{\theta}}_{\phi,n}^{(m)}-\eta) \stackrel{d}{\rightarrow} 
			N\left(0,\ {I(\eta)}^{-1}\right),\text{ as }n\to\infty. 
		\end{align*}
	\end{itemize}
	
\end{thm}
Theorem 1 can be proved using the Cram\'er-Wald device, and arguments similar to those used in proving Theorem 2 of Ekstr\"om et al. (2020). Details of the proof are provided in the supplementary material. 

We have the following result concerning the asymptotic distribution of the test statistic $\tilde{T}_{\phi,n}^{(m)}(\theta_0)$, under the null hypothesis and under a sequence of local alternatives around $\theta_0$.
\begin{thm}\label{thm2} 
	Let $m$ be a fixed positive integer. Assume that  $\hat{\theta}_{\phi,n}^{(m)}= \arg \inf_{\theta\in\Theta}S_{\phi,n}$ is a consistent estimator of the true parameter $\eta\in\Theta\subseteq\mathbb{R}^p$. Suppose that assumptions of Theorem 1 hold. Then,
	\begin{itemize}
		\item[(a) ] under $H_0:\eta=\theta_0$, $\tilde{T}_{\phi,n}^{(m)}(\theta_0)\stackrel{d}{\rightarrow}\chi^2_p$, as $n\to\infty$; and
		\item[(b) ] under ${\eta}_{n}=\theta_0+\Delta n^{-1/2}$, where $\Delta\in\mathbb{R}^p$ is a fixed vector (independent of $n$), $\tilde{T}_{\phi,n}^{(m)}(\theta_0)\stackrel{d}{\rightarrow}\chi^2_p \left(\sigma_{\phi,m}^{-2} \Delta^t I(\theta_0)\Delta\right)$, as $n\to\infty$.
	\end{itemize}
\end{thm}
The special case of Theorem \ref{thm2}, when $m=1$ was proved by Ekstr\"om (2013).
Thus, for testing (\ref{simpleh0}) based on the statistic $\tilde{T}_{\phi,n}^{(m)}(\theta_0)$, we may take rejection region as $\omega=\{\tilde{T}_{\phi,n}^{(m)}(\theta_0)\geq c_\alpha\}$, where $\alpha$ is the size of the test and the critical value $c_\alpha$ is the $(1-\alpha)^{th}$ quantile of the ${\chi_p^2}$ distribution. The power of the test for local alternatives ${\eta}_{n}=\theta_0+\Delta n^{-1/2}$, around $\theta_0$, is given by 
\begin{align*}
	\beta^*({\theta}_0)=\int_{c_\alpha}^{\infty}\ dF_{\chi_p^2\left(\sigma_{\phi,m}^{-2} \Delta^t I(\theta_0)\Delta\right)}(y),
\end{align*}
where $F_{\chi_p^2\left(\sigma_{\phi,m}^{-2} \Delta^t I(\theta_0)\Delta\right)}$ denotes the distribution function of ${\chi_p^2\left(\sigma_{\phi,m}^{-2} \Delta^t I(\theta_0)\Delta\right)}$.
To maximise the power function, it is evident that $\phi$ should be such that $\sigma_{\phi,m}^2$ is minimum. 
Ekstr\"om et al. (2020) found that $\sigma_{\phi,m}^2\geq 1$ for any $\phi$ and $m$. Further, they observed that for any finite $m$, $\sigma_{\phi,m}^2$ attains the minimum value if and only if $\phi(x)=a\log(x)+b\ x+c,\ x>0$, for some real constants $a~(<0)$, $b$ and $c$. Thus, for finite $m$, within the class of tests based on $\tilde{T}_{\phi,n}^{(m)}(\theta_0)$, the test corresponding to $\phi(x)=-\log(x),\ x>0,$ is asymptotically locally most powerful. For $m\to\infty$, we have the following result.
\begin{thm}\label{thm3} Suppose that $m\to\infty$ such that $m=o(n)$, and $\lim_{m\to\infty}\sigma_{\phi,m}^2=1$.	Assume that  $\hat{\theta}_{\phi,n}^{(m)}= \arg \inf_{\theta\in\Theta}S_{\phi,n}^{(m)}(\theta)$ is a consistent estimator of the true parameter $\eta\in\Theta\subseteq\mathbb{R}^p$. Further suppose that assumptions of Theorem \ref{thm1} hold. Then,
	\begin{itemize}
		\item[(a) ] under $H_0:\eta=\theta_0$, $\tilde{T}_{\phi,n}^{(m)}(\theta_0)\stackrel{d}{\rightarrow}\chi^2_p$, as $n\to\infty$; and
		\item[(b) ] under ${\eta}_{n}=\theta_0+\Delta n^{-1/2}$, $\tilde{T}_{\phi,n}^{(m)}(\theta_0)\stackrel{d}{\rightarrow}\chi^2_p \left(\Delta^t I(\theta_0)\Delta\right)$, as $n\to\infty$.
	\end{itemize}
\end{thm}
The functions $\phi$, defined in (\ref{pow_fun}), satisfy the assumption $\lim_{m\to\infty}\sigma_{\phi,m}^2=1$ (cf. Ekstr\"om et al. 2020). 
Theorem \ref{thm3} suggests that if $m=o(n)$, then all the tests based on $\phi$, with $\lim_{m\to\infty}\sigma_{\phi,m}^2=1$, have the same asymptotic power. 

Observe that we can use two different convex functions $\phi_1$ and $\phi_2$ in (\ref{teststt}), for the purpose of estimating and testing  {(cf. Ekstr\"om 2013)}. By doing so, we obtain a larger class of tests based on the following test statistic
\begin{align}\label{teststt1}
	\tilde{T}_{\phi_1,\phi_2,n}^{(m)}(\theta_0)= \frac{{2n}\left(S_{\phi_1,n}^{(m)}(\theta_0) -S_{\phi_1,n}^{(m)}(\hat{\theta}_{\phi_2,n}^{(m)})\right)}{E\left({\zeta}_{m}^{2} \phi_1 ''({\zeta}_{m} )\right)\sigma_{\phi_2,m}^2}.
\end{align}
By using $\phi_2(x)={-\log}(x)$ in (\ref{teststt1}), we get an asymptotically locally efficient testing procedure based on sample spacings for fixed $m$.
When $m=o(n)$, by choosing any $\phi_2$ such that $\lim_{m\to\infty}\sigma_{\phi_2,m}^2=1$, we get an asymptotically locally efficient testing procedure based on sample spacings. Define
\begin{align*}
	{\tilde{T}}^{(m)*}_{\phi,n}(\theta_0):=\frac{2n\left(S_{\phi,n}^{(m)}(\theta_0) -S_{\phi,n}^{(m)}(\hat{\theta}_{-\log,n}^{(m)})\right)} {E\left({\zeta}_{m}^{2} \phi ''({\zeta}_{m} )\right)}.
\end{align*}
It is possible to extend Theorems \ref{thm2} and \ref{thm3} to the statistic in (\ref{teststt1}), as stated in following results:
\begin{cor}\label{cor1}
	Suppose that $m$ is finite and that assumptions of Theorem \ref{thm2} hold. Then,
	\begin{itemize}
		\item[(a) ] under $H_0$, ${\tilde{T}}^{(m)*}_{\phi,n}(\theta_0)\stackrel{d}{\rightarrow}\chi^2_p$, as $n\to\infty$; and
		\item[(b) ] under ${\eta}_{n}=\theta_0+\Delta n^{-1/2}$, ${\tilde{T}}^{(m)*}_{\phi,n}(\theta_0)\stackrel{d}{\rightarrow}\chi^2_p \left( \Delta^t I(\theta_0)\Delta\right)$, as $n\to\infty$.
	\end{itemize}
\end{cor}
\begin{cor}\label{cor2}
	Let $m\to\infty$ such that $m=o(n)$, and $\lim_{m\to\infty}\sigma_{\phi_2,m}^2=1$. Further, suppose that the assumptions of Theorem \ref{thm3} hold. Then
	\begin{itemize}
		\item[(a) ] under $H_0$, $\tilde{T}_{\phi_1,\phi_2,n}^{(m)}(\theta_0)\stackrel{d}{\rightarrow}\chi^2_p$, as $n\to\infty$; and
		\item[(b) ] under ${\eta}_{n}=\theta_0+\Delta n^{-1/2}$, $\tilde{T}_{\phi_1,\phi_2,n}^{(m)}(\theta_0)\stackrel{d}{\rightarrow}\chi^2_p \left( \Delta^t I({\theta_0})\Delta\right)$, as $n\to\infty$.  
	\end{itemize}
\end{cor}
Using Corollaries \ref{cor1} and \ref{cor2}, it follows that tests based on ${\tilde{T}}^{(m)*}_{\phi,n}$ and $\tilde{T}_{\phi_1,\phi_2,n}^{(m)}(\theta_0)$, such that $\lim_{m\to\infty}\sigma_{\phi_2,m}^2=1$, are asymptotically equivalent to the LRT (Sen and Singer 1994). 

Here, $\omega_\alpha=\{\tilde{T}_{\phi,n}^{(m)}(\theta_0)\geq c_\alpha\}$ is a critical region of size $\alpha$, where $c_\alpha$ is the $(1-\alpha)^{th}$ quantile of ${\chi_p^2}$ distribution. We can also adopt the p-value approach for parametric tests. It can be easily seen that, for testing (\ref{simpleh0}), tests based on $\tilde{T}_{\phi,n}^{(m)}(\theta_0)$ form nested critical regions, i.e., $\omega_\alpha\subset\omega_{\alpha'}$ if $\alpha<\alpha'$. In case of nested critical regions, it would be useful to determine not only whether the verdict is to reject the null hypothesis at the given level of significance, but also to determine the least level of significance (i.e., $\hat{p}(X)=\inf\{\alpha:X\in\omega_\alpha\}$)  at which the null hypothesis would be rejected for the available observation set (cf. Lehman and Romano 2005, p~63). With the p-value approach, one can also determine whether the decision to reject the null hypothesis is made convincingly, or it is just a borderline decision.
\subsection{Composite Null Hypothesis}
The testing procedure discussed in the previous section can also be extended for testing composite hypotheses. For testing composite hypotheses it is pragmatic to present testing problem as  
\begin{align*}
	H_0:\ h(\eta)=0\text{ against }H_A:\ h(\eta)\neq0,
\end{align*}
where $\Theta\subseteq\mathbb{R}^p$ and $h=(h_1,h_2,\ldots,h_r):\mathbb{R}^p\to\mathbb{R}^r$ is a vector-valued function such that for every $\theta\in\Theta$, the $p\times r$ matrix, $H(\theta)=((\frac{\partial}{\partial\theta_i}h_j(\theta)))$ exists, $H(\theta)$ is continuous in $\theta$ and $rank\big(H(\theta)\big)=r$. For example, if $p=3$ and $\eta=(\eta_1,\eta_2,\eta_3)^t$, then $H_0:\eta_1=\eta_2$ corresponds to $h(\eta)=\eta_1-\eta_2$.
Under the above setup, an extension of the test statistic (\ref{teststt}) is 
\begin{align*}
	{T}_{\phi,n}^{(m)}(h)=\frac{{2n}\left(\inf_{\theta:h(\theta)=0}S_{\phi,n}^{(m)}(\theta)-\inf_{\theta\in\Theta}S_{\phi,n}^{(m)}(\theta)\right)} {E\left({\zeta}_{m}^{2} \phi ''({\zeta}_{m} )\right)\sigma_{\phi,m}^2}.
\end{align*}
Significant large positive value of the test statistic ${T}_{\phi,n}^{(m)}(h)$ would result into the rejection of the null hypothesis. We have the following result concerning asymptotic behaviour of ${T}_{\phi,n}^{(m)}(h)$, under the null hypothesis.
\begin{thm}\label{thm4} Let $m=o(n)$ (finite or $m\to\infty)$ and suppose that assumptions of  Theorem \ref{thm1} hold. Assume that  $\hat{\theta}_{\phi,n}^{(m)}= \arg \inf_{\theta\in\Theta}S_{\phi,n}$ is a consistent estimator of the true parameter $\eta\in\Theta\subseteq\mathbb{R}^p$. Then, under $H_0$, 
	\begin{align*}
		{T}_{\phi,n}^{(m)}(h)\xrightarrow{d}\chi^2_{r},\text{ as } n\to\infty.
	\end{align*}
\end{thm}
The above result gives another similarity between spacings based tests and the LRT (Sen and Singer 1994). Observe that the test procedure based on $\hat{T}_{\phi,n}^{(m)}(h)$ can be used for testing hypotheses concerning the regression coefficients of a regression model when errors follow a mixture, or a heavy-tailed distribution and the LRT is not applicable.  

This idea has also been extended to nonparametric setting (Del Pino 1979). A related discussion can be found in Ekstr\"om (2013).

\section{Tests for Multivariate Distributions}
The results of the previous section can be extended to multivariate distributions for $m=1$. The definition of spacings for multivariate observations can be defined in terms of nearest neighbor balls. Let $X_1,X_2,\ldots,X_{n}$ be independent and identically distributed $d$-dimensional random vectors from an absolutely continuous distribution function $F_\eta$, $\eta\in \Theta\subseteq \mathbb{R}^p$. 
The nearest neighbour distance to $X_i$ is defined as 
\begin{align}
	R_n(i):= \min_{i\neq j}||X_i-X_j||,\ i=1,2,\ldots,n,\nonumber
\end{align}
where $||\cdot||$ is some distance measure on $\mathbb{R}^d$. Let $B(x,r)=\{y\in\mathbb{R}^d:||x-y||\leq r\}$ denotes the closed ball with centre $x\in\mathbb{R}^d$ and radius $r$ $(>0)$.  
Denote the nearest neighbour to $X_i$ by $X_i(nn)$ and the nearest neighbour ball of $X_i$ by $B_n(X_i):= \{y\in\mathbb{R}^d:||X_i-y||\leq R_n(i)\}$. Let $f_\theta$ and $\mathbb{P}_\theta$, respectively, denote the density function and the probability measure corresponding to the distribution $F_\theta$, $\theta\in\Theta$. Define a random variable $\xi_{i,n}(\theta)$ for $\theta\in\Theta$, as
\begin{align}
	\xi_{i,n}(\theta)= n\mathbb{P}_\theta(B_n(X_i))=n\int_{B_n(X_i)}\ d~\mathbb{P}_\theta(y),~  i=1,2,\ldots,n\nonumber.
\end{align}
Ranneby et al. (2005) extended the idea of maximum spacing estimator to multivariate observations and defined the multivariate maximum spacing estimator as $\hat{\theta}_n=\arg \sup_{\theta\in\Theta} \frac{1}{n}\sum_{i=1}^{n}\log(\xi_{i,n}(\theta))$. Kuljus and Ranneby (2015) extended this idea to any strictly convex function $\phi:(0,\infty)\to(0,\infty)$ such that $\phi$ has the minima at $x=1$.  They defined the generalised spacing function $S_{\phi,n}(\theta)$ and the generalised spacing estimator $\hat{\theta}_{\phi,n}$ as follows:
\begin{align*}
	S_{\phi,n}(\theta)=\frac{1}{n}\sum_{i=1}^{n}\phi(\xi_{i,n}(\theta)),\theta\in\Theta \text{ and }  \hat{{\theta}}_{\phi,n}=\arg \min_{\theta\in\Theta} S_{\phi,n}(\theta).
\end{align*}
If the minimiser does not exist, then the estimator can be suitably modified on the lines suggested by Kuljus and Ranneby (2015). For simplicity, we assume that the minimiser exists. 

Consider following notations:
\begin{align*}
	&q(x)=x\phi'(x),\ \sigma_q^2=q^2(0)+\int_{0}^{\infty}\int_{0}^{\infty} k(s,t)dq(s)dq(t)+2q(0)\int_{0}^{\infty} k(0,t)dq(t),\\
	&\text{and }k(s,t)=e^{-t}-te^{-s-t}+e^{-s-t}\int_{W(s,t)}(e^{\beta(s,t,x)}-1)dx,\ 0\leq s\leq t\leq \infty.
\end{align*}
Here,
\begin{align*}
	W(s,t)=\{x\in\mathbb{R}^d: r_1\leq ||x||\leq r_1+r_2\}\text{ and } \beta(s,t,x)=\int_{B(0,r_1)\cap B(x,r_2)}dz,
\end{align*}
with $t$ and $s$ denoting volumes of the balls $B(0,r_1)$ and $B(0,r_2)$, respectively.

 {Let $\theta^{(n)}_0=~\arg\min_{\theta\in\Theta}\mathbb{E}[\phi(\xi_{1,n}(\theta))]$}. Then, under some general assumptions, Kuljus and Ranneby (2020) proved the asymptotic normality of $ \hat{\theta}_{\phi,n}$, i.e.,
\begin{align*}
	\sqrt{n}(\hat{\theta}_{\phi,n}-\theta^{(n)}_0)\stackrel{d}{\rightarrow}N\left(0,\frac{\sigma_q^2}{b_\phi^2}I(\eta)^{-1} \right),\text{ as }n\to\infty,
\end{align*}
where $b_\phi=\mathbb{E}(Z_1^2\phi''(Z_1))$, $\sigma_q^2$ is defined above, and $I(\eta)$ is the Fisher information matrix.

For a specified $\theta_0\in\Theta$, we wish to test the null hypothesis $H_0:\eta=\theta_0$ against the alternative $H_A:\eta\neq\theta_0$.
Define, 
\begin{align}
	T_{\phi,n}(\theta_0) &=\frac{1}{n}\sum_{i=1}^{n}\phi(\xi_{i,n}(\theta_0))-\inf_{\theta\in\Theta}\frac{1}{n}\sum_{i=1}^{n}\phi(\xi_{i,n}(\theta))\nonumber\\
	&=\frac{1}{n}\sum_{i=1}^{n}\phi(\xi_{i,n}(\theta_0))-\frac{1}{n}\sum_{i=1}^{n}\phi(\xi_{i,n}(\hat{\theta}_{\phi,n}))\nonumber\\
	\text{and }\tilde{T}_{\phi,n}(\theta_0)&= \frac{2nb_\phi}{\sigma_q^2}T_{\phi,n}.\label{multstt}
\end{align}
For testing $H_0: \eta=\theta_0$ against $H_A:\eta\neq\theta_0$, on the lines of earlier discussion, we will reject $H_0$ for large values of the test statistic $\tilde{T}_{\phi,n}(\theta_0)$. The following theorem provides the asymptotic distribution of $\tilde{T}_{\phi,n}(\theta_0)$, under the null hypothesis and a sequence of local alternatives.
\begin{thm}\label{thm6}
	Assume that $\hat{\theta}_{\phi,n}=\arg \min_{\theta\in\Theta} S_{\phi,n}(\theta)$ is a consistent estimator of true $\eta\in\Theta\subseteq\mathbb{R}^p$. Suppose that $\sqrt{n}(\theta_{0}^{(n)}-\theta_{0})\xrightarrow{p}0$ as $n\to\infty$, and that assumptions of Theorem 1 of Kuljus and Ranneby (2020) hold. Then,
	\begin{itemize}
		\item[(a) ] under $H_0$, $\tilde{T}_{\phi,n}(\theta_0)\stackrel{d}{\rightarrow}\chi^2_p$ as $n\to\infty$; and
		\item[(b) ] under ${\eta}_{n}=\theta_0+\Delta n^{-1/2}$, $\tilde{T}_{\phi,n}(\theta_0)\stackrel{d}{\rightarrow}\chi^2_p \left(\dfrac{b_\phi^2}{\sigma_q^2} \Delta^t I(\theta_0)\Delta\right)$ as $n\to\infty$.
	\end{itemize}
\end{thm}
The assumption $\sqrt{n}(\theta_{0}^{(n)}-\theta_{0})\xrightarrow{p}0$ as $n\to\infty$, in Theorem \ref{thm6} may be sometimes difficult to verify. But, it holds in some commonly encountered situations, e.g., in the $d$-variate normal distribution $N_d(\theta_0,\Xi_0)$, $\theta_{0}^{(n)}=\theta_{0}$ (cf. Kuljus and Ranneby 2020).

Using Theorem \ref{thm6}, an asymptotic test procedure of level $\alpha$ would be exactly the same as in the univariate case. Also, the asymptotic power of this test, for local alternatives, is again the same as the univariate case. 
This procedure can be extended for testing composite hypotheses concerning absolutely continuous  multivariate distributions. 
Practically, a composite hypothesis testing problem can be written as
\begin{align*}
	H_0:\ h(\eta)=0,\text{ against, }H_A:\ h(\eta)\neq0,
\end{align*}
where $h=(h_1,h_2,\ldots,h_r):\mathbb{R}^p\to\mathbb{R}^r$ is a vector valued function such that the $p\times r$ matrix, $H(\theta)=((\frac{\partial}{\partial\theta_i}h_j(\theta)))$ exists and is continuous in $\theta$, and $rank\big(H(\theta)\big)=r$.
An extension of the test is based on the statistic
\begin{align*}
	{T}_{\phi,n}(h)=\frac{2nb_\phi}{\sigma_q^2} \left(\ \inf_{\theta:h(\theta)=0}\frac{1}{n}\sum_{i=1}^{n}\phi(\xi_{i,n}(\theta))-\inf_{\theta\in\Theta}\frac{1}{n}\sum_{i=1}^{n}\phi(\xi_{i,n}(\theta)) \right),
\end{align*}
where significantly large positive values of the test statistic would result into rejection of the null hypothesis. 
\begin{thm}\label{thm7}
	Assume that $\hat{\theta}_{\phi,n}=\arg \inf_{\theta\in\Theta} S_{\phi,n}(\theta)$ be a consistent estimator of true $\eta\in\Theta\subseteq\mathbb{R}^p$, and $\sqrt{n}(\theta_{0}^{(n)}-\eta)\xrightarrow{p}0$ as $n\to\infty$. Suppose that assumptions of Theorem 1 of Kuljus and Ranneby (2020) hold. Then, under $H_0$, we have
	\begin{align*}
		{T}_{\phi,n}(h)\xrightarrow{d}\chi^2_{r}\ as\ n\to\infty.
	\end{align*}
\end{thm}
Similar to the univariate case, we can use two different convex functions $\phi_1$ and $\phi_2$, one for the estimation and other for the testing $H_0$, and obtain results similar to ones stated in Theorems \ref{thm6} and \ref{thm7}.

Let $\rho$ be a semi-metric, i.e., it satisfies the following conditions: $\rho(s,t)\geq0$, $\forall$ $s$ and $t$; $\rho(s,t)=0$, iff $s=t$; and $\rho(s,t)=\rho(t,s)$, $\forall$ $s$ and $t$. Let $\rho(s,\cdot)$ be thrice continuously differentiable and that $\tau_{\phi,m}=\frac{\partial}{\partial s}\rho(s,\mathbb{E}(\zeta_m))|_{s=\mathbb{E}(\zeta_m)}$ exists and is non-zero. Then, all the test statistics (univariate and multivariate) discussed above can be modified by replacing the Euclidean distance by the semi-metric, and versions of theorems and corollaries stated above can be obtained (cf. Ekstr\"om 2013).

Yu and Ekstr\"om (2000) studied spacings based on identically distributed univariate observations (not necessarily independent), satisfying certain mixing conditions. They proposed an estimator based on such spacings. The proposed estimator consistently estimates the entropy and the argument maximizing the estimate of the entropy is a consistent estimate of unknown parameter. Asymptotic distribution of such an estimators and parametric tests based on such spacings are still unresolved. 
\section{Numerical Study}
In this section we will evaluate finite sample performances of the proposed tests and compare their empirical powers with those of the LRT. One of the essential properties of a good test is that its Type-I error rate should be close to the level of the test. If the Type-I error rate is away from the level then the conclusion made may be spurious.  Obviously, both, the type-I error rate and the power of the test, depend on both $n$ and $m$. In the following numerical study we have taken $10,000$ Monte-Carlo samples to estimate type-I error rates and powers of various tests.
\subsection{Univariate Case}
\begin{figure}[h]
	\centering
	\includegraphics[width=7cm,height=6.5cm]{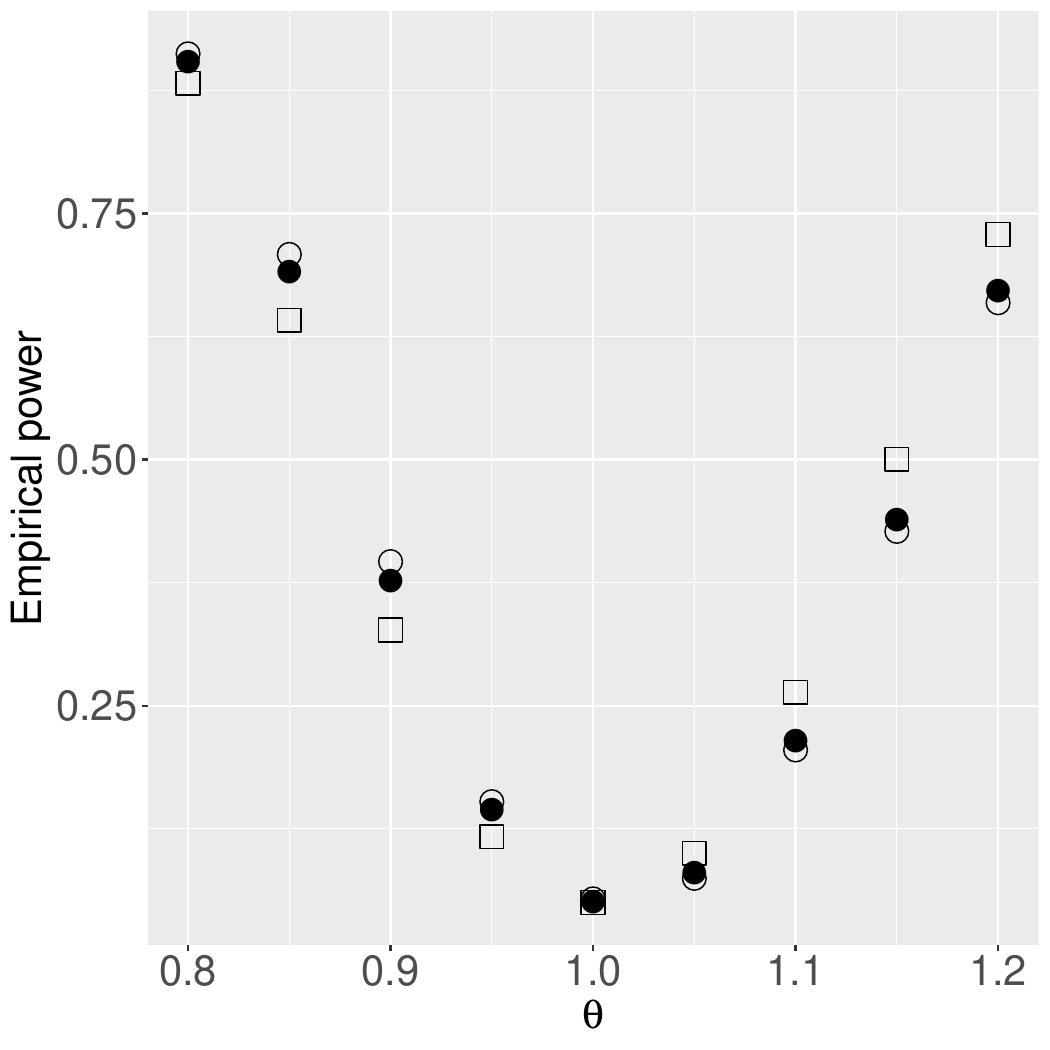}
	\caption{The empirical powers of $\tilde{T}_{-\log,n}^{(1)}(1)$ (hollow circles), $\tilde{T}_{-\log,n}^{(m_{opt}(n))}(1)$ (solid circles) and LRT (hollow squares), for the alternative $\eta=\theta$ and $n=200$.}\label{fig1}
\end{figure}

For a given sample size $n$, let $m_{opt}(n)$ denote  the value of $m$ for which the Type-I error rate is closest to the level of the test. Note that $m_{opt}(n)$ may also depend on the hypothesized model. To obtain $m_{opt}(n)$, we propose the following algorithm based on the hypothesized model and data.
\begin{itemize}
	\item[Step-I]  Compute $\theta_{\phi,n}^{(1)}$ (or, $\theta_{\phi,n}^{(m)}$ for any fixed $m$) using the given data and the hypothesized model.
	\item[Step-II] For some large positive integer $B$, at the $b^{th}$ occasion ($b=1,2,\ldots,B$), draw a bootstrap sample $x_{b1}^*,x_{b2}^*,\ldots,x_{bn}^*$ from the population with distribution $F_{\hat{\theta}_{\phi,n}^{(1)}}$. Let $\mathcal{M}$ be a subset of suitably chosen natural numbers, containing competing choices of $m$. For each $m\in\mathcal{M}$, compute the Type-I error rate for the test, based on the $B$ bootstrap samples.
	\item[Step-III] Take $m_{opt}(n)=\big\{m\in\mathcal{M}:$ Type-I error rate is closest to the level for the test based on $\tilde{T}_{\phi,n}^{(m)}\big\}$. 
\end{itemize}
We have taken $\mathcal{M}$ to be the set of all natural numbers $\leq\frac{n}{10}$. If $m_{opt}(n)$ has multiple values, then we take minimum of those $m_{opt}(n)$.\\~\\
\textit{Example 1.} Let $X_1,\ldots,X_{n}$ be an available random sample from an exponential distribution with failure rate parameter $\eta\in (0,\infty)=\Theta$. We want to test $H_0:\eta=1$ against the alternative $H_A:\ \eta\neq 1$. Under $H_0$ and $m=o(n)$, $\tilde{T}_{\phi,n}^{(m)}(1)\stackrel{d}{\rightarrow}\chi^2_1$ as $n\to\infty$. We take $\alpha=0.05$ as the level of significance. Then, $H_0$ is rejected if $\tilde{T}_{\phi,n}^{(m)}>\chi^2_{1,0.05}\simeq 3.841$. Ekstr\"om (2013) has investigated the performance of this test for $n=100$, $m=1$ and $\phi(x)=-\log(x)$. We have studied performance of the test procedure based on $\tilde{T}_{\phi,n}^{(m)}(1)$, for different values of $n$ and $m$, and fixed alternatives $\eta\in\{0.80,0.85,0.90,0.95,1.00,1.05,1.15,1.10,1.15,1.20\}$. We found that $m_{opt}(n)$ increases as $n$ increases, and, for $n\leq 1200$, $1\leq m_{opt}(n)\leq 5$. Figure \ref{fig1} gives empirical powers of the test for $n=200$ and $\phi(x)=-\log(x)$. Here, we have $m_{opt}(n)=3$.
\FloatBarrier
When the failure rate is greater than $1$, the LRT performs better than tests based on $\tilde{T}_{-\log,n}^{(m)}(1)$, but for failure rate less than $1$, tests based on $\tilde{T}_{-\log,n}^{(m)}(1)$ outperform the LRT.  {For some other choices of $\phi$ (e.g., $\phi_\gamma$ with $\gamma=-0.9,-0.5$), tests based on $\tilde{T}_{\phi,n}^{(m)}(1)$ perform slightly inferior.}

\begin{figure}[h]
	\centering
	\begin{subfigure}[t]{0.48\textwidth}
		\includegraphics[width=6cm,height=6cm]{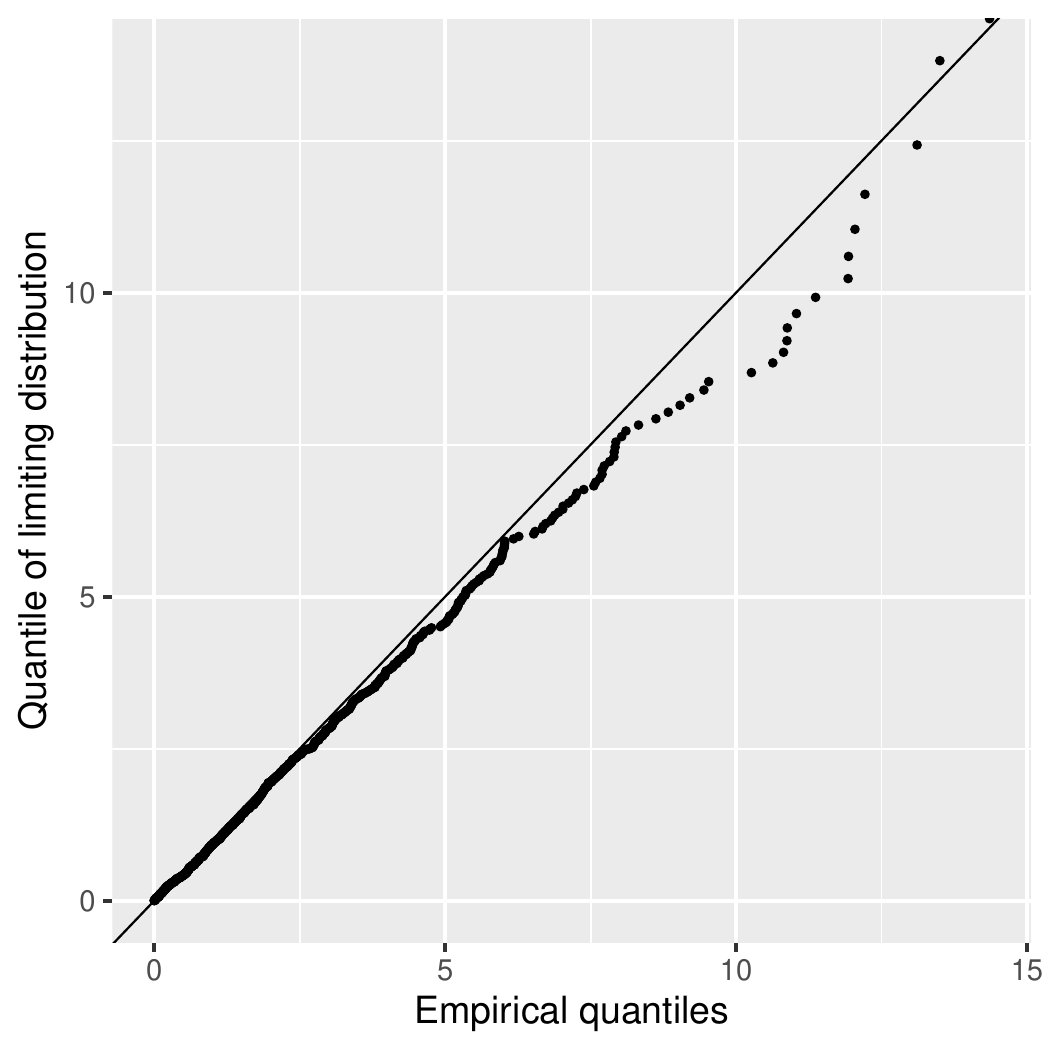}
		\caption*{$m=1$}
	\end{subfigure}
	\begin{subfigure}[t]{0.48\textwidth}
		\includegraphics[width=6cm,height=6cm]{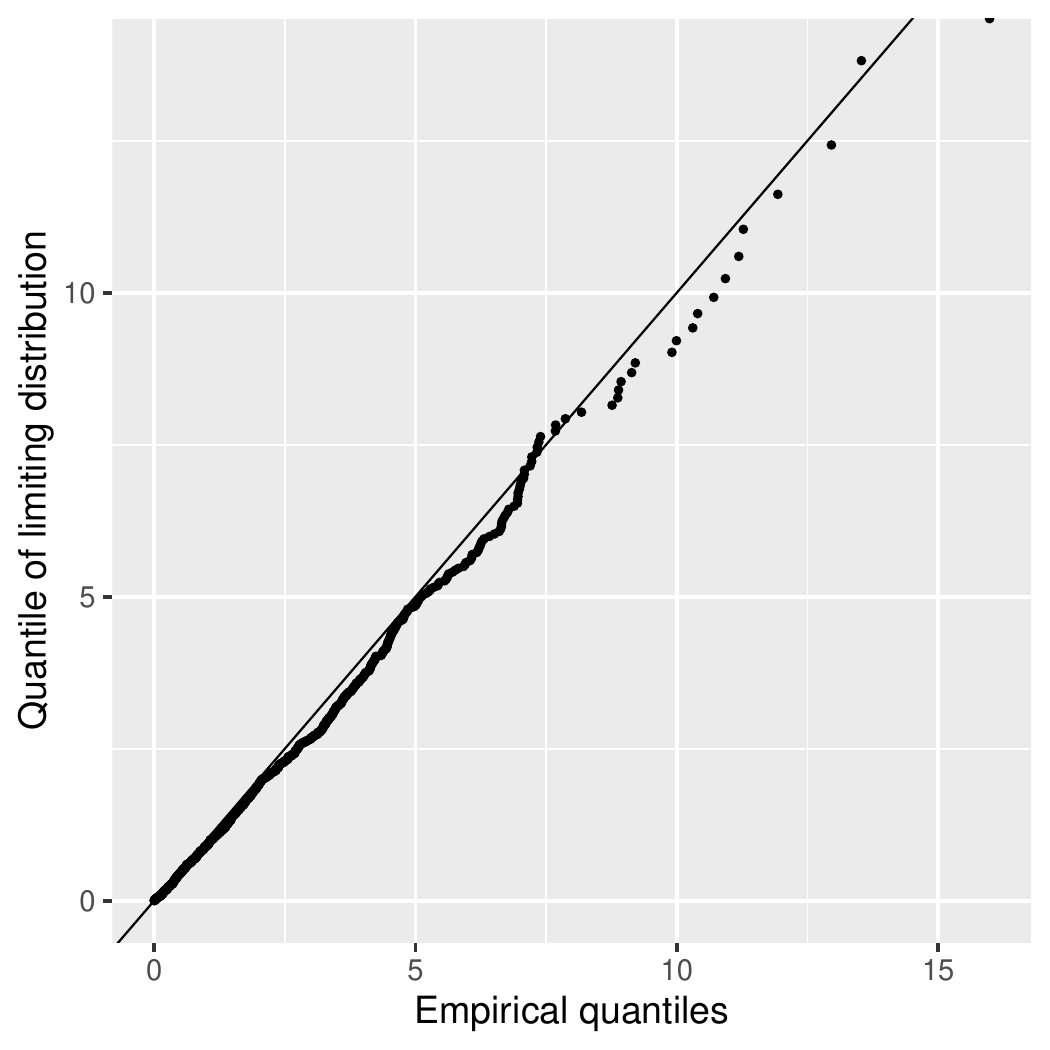}
		\caption*{$m=2$}
	\end{subfigure}
	\caption{Q-Q plot of 1000 replicates of $\tilde{T}_{-\log,n}^{(m)}(\theta_0)$,  $n=225$ under $H_0$. Empirical distribution quantiles on the horizontal axis and the limiting $\chi^2_2$ distribution quantiles on the vertical axis.}\label{fig2}
\end{figure}

\begin{figure}[h]
	\centering
	\begin{subfigure}[t]{0.48\textwidth}
		\includegraphics[width=6cm,height=6cm]{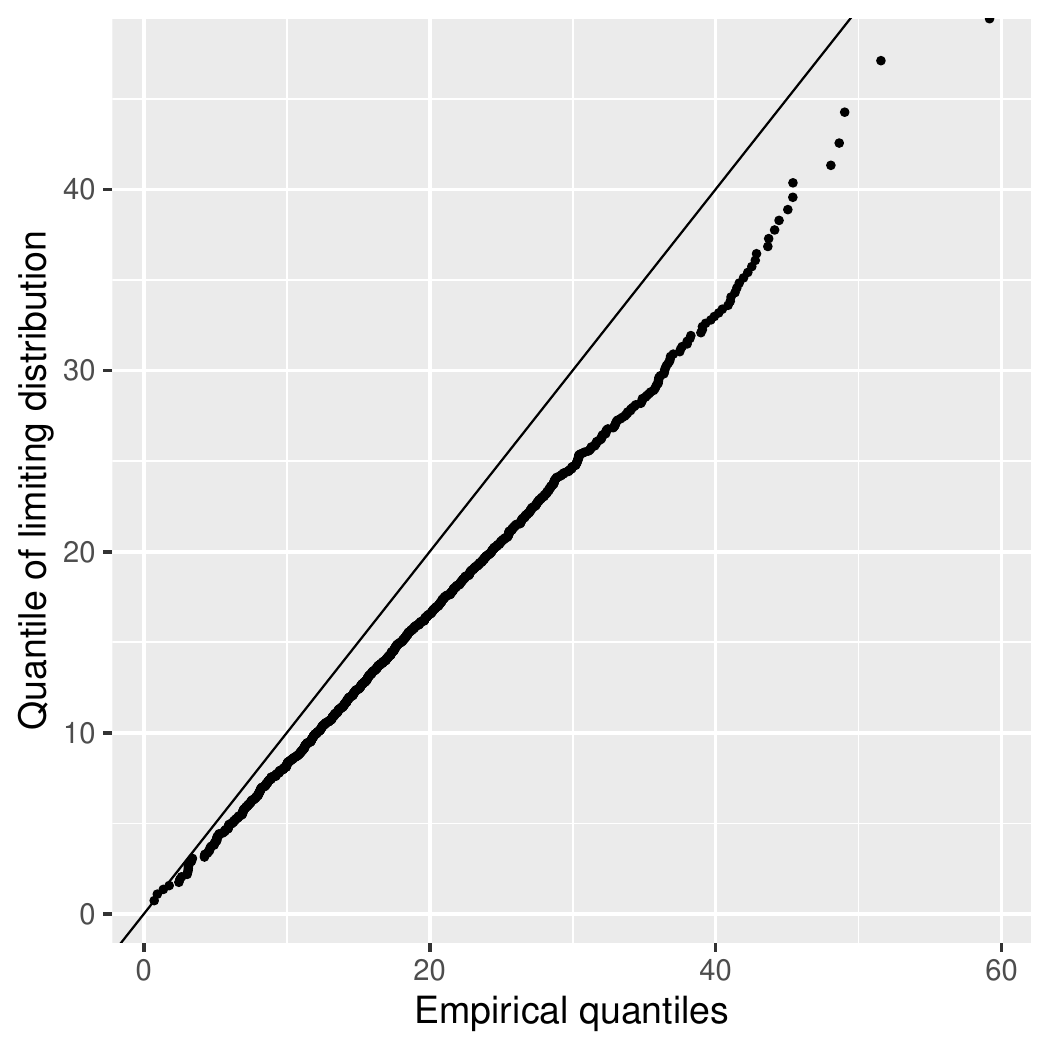}
		\caption*{$m=1$}
	\end{subfigure}
	\begin{subfigure}[t]{0.48\textwidth}
		\includegraphics[width=6cm,height=6cm]{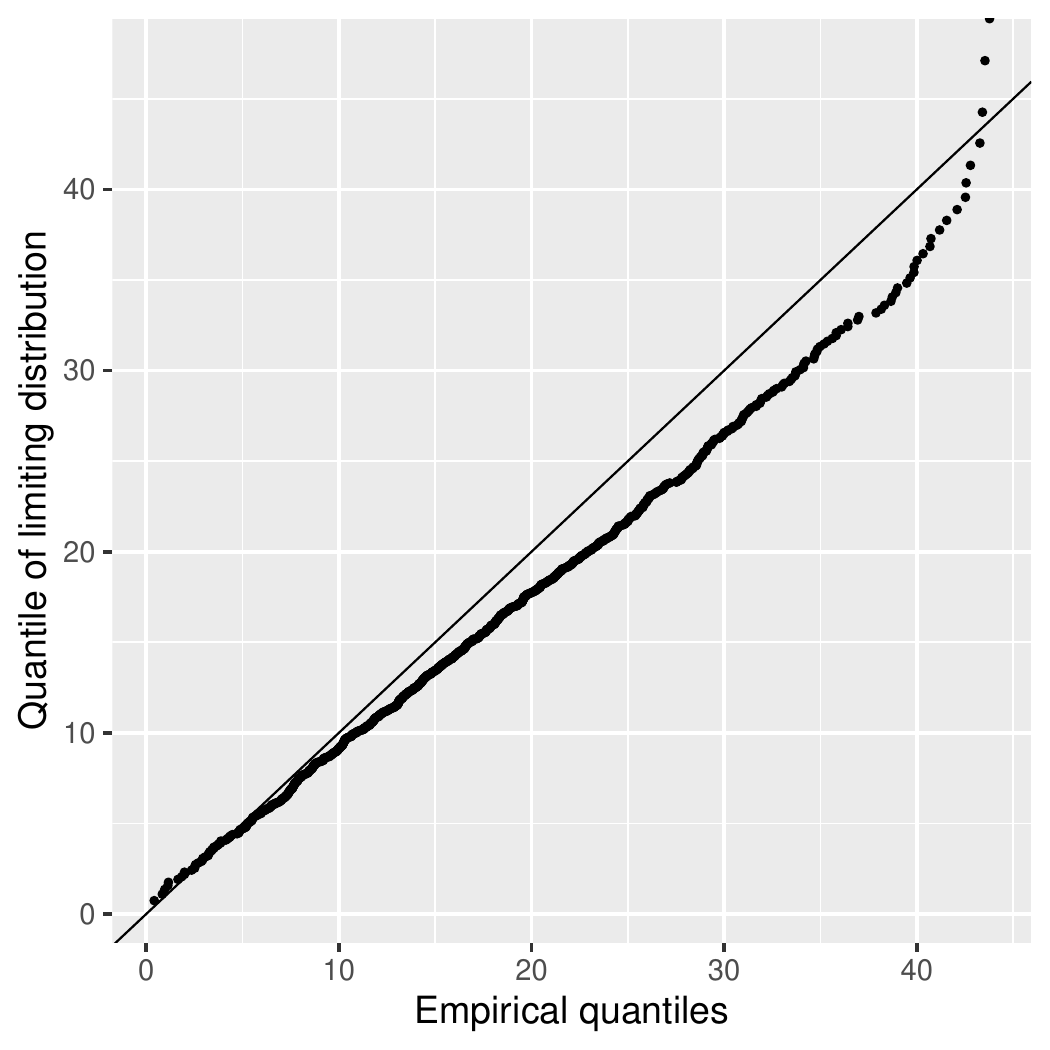}
		\caption*{$m=2$}
	\end{subfigure}
	\caption{Q-Q plot of 1000 replicates of $\tilde{T}_{-\log,n}^{(m)}(\theta_0)$,  $n=225$ under the alternative $\eta_n=\theta_0+\Delta n^{-1/2}$ with $\Delta=(3,3)$. Empirical distribution quantiles on the horizontal axis and the limiting $\chi^2_2(\Delta^tI(\theta_0)\Delta)$ distribution quantiles on the vertical axis.}\label{fig3}
\end{figure}

If assumptions (A1)-(A6) are not satisfied, then the limiting distribution of $\tilde{T}_{\phi,n}^{(m)}$ may not be chi-square (cf. Ekstr\"om 2013). Sometimes, due to unboundedness of the likelihood function, the LRT can not be used but tests based on $\tilde{T}_{\phi,n}^{(m)}(\theta_0)$ may be applicable, as illustrated in the next example.\\~\\
\textit{Example 2.} Consider a population with a distribution function,
\begin{align*}
	F_\eta(x)=\frac{1}{2}\Phi(x-\mu)+\frac{1}{2}\Phi\left(\frac{x-\mu}{\sigma}\right),\ -\infty <x< \infty,
\end{align*}
where $\eta=(\mu,\sigma)\in\mathbb{R}\times(0,\infty)$. Let $X_1,\ldots,X_{n}$ be a random sample from this distribution. Suppose that we want to test $H_0:\eta=(0,1)$ against $H_A:\eta\neq (0,1)$. For this problem, conditions required for the proposed tests are satisfied. Ekstr\"om (2013) observed that the LRT breaks down in this situation. Fix $n=225$. Figure \ref{fig2} shows that the distribution of $\tilde{T}_{-\log,n}^{(m)}(0,1)$ is better approximated by the limiting chi-square distribution for $m=m_{opt}=2$ than for $m=1$. For local alternative $\eta_n=(0,1)+\Delta n^{-1/2}$, we take $\Delta=(3,3)$. Then, for $n=225$, $\eta_n=(0.2,1.2)$. Under this alternative, Figure \ref{fig3} shows that  the distribution of $\tilde{T}_{-\log,n}^{(m)}(0,1)$, for $m=m_{opt}(n)=2$, is closer to the limiting non-central chi-square distribution than for the case of $m=1$. For other choices of  $\phi$ (as in Example 1) similar observations were made.

We observed that $m_{opt}(n)$ depends not only on $n$ but also on the underlying model. 	 For the models considered in our study, as $n$ increases the type-I error rates remain close to the level upto certain large values of $m$. For sample size $n\geq50$, our proposed tests, with $m=m_{opt}(n)$, perform similar to the likelihood ratio test. We also observed that the test corresponding to $\phi(x)=-\log(x)$ performs better than tests corresponding to other $\phi$s in (\ref{pow_fun}).
\subsection{Multivariate Case}
\begin{figure}[h]
	\centering
	\begin{subfigure}[t]{0.48\textwidth}
		\includegraphics[width=6cm,height=6cm]{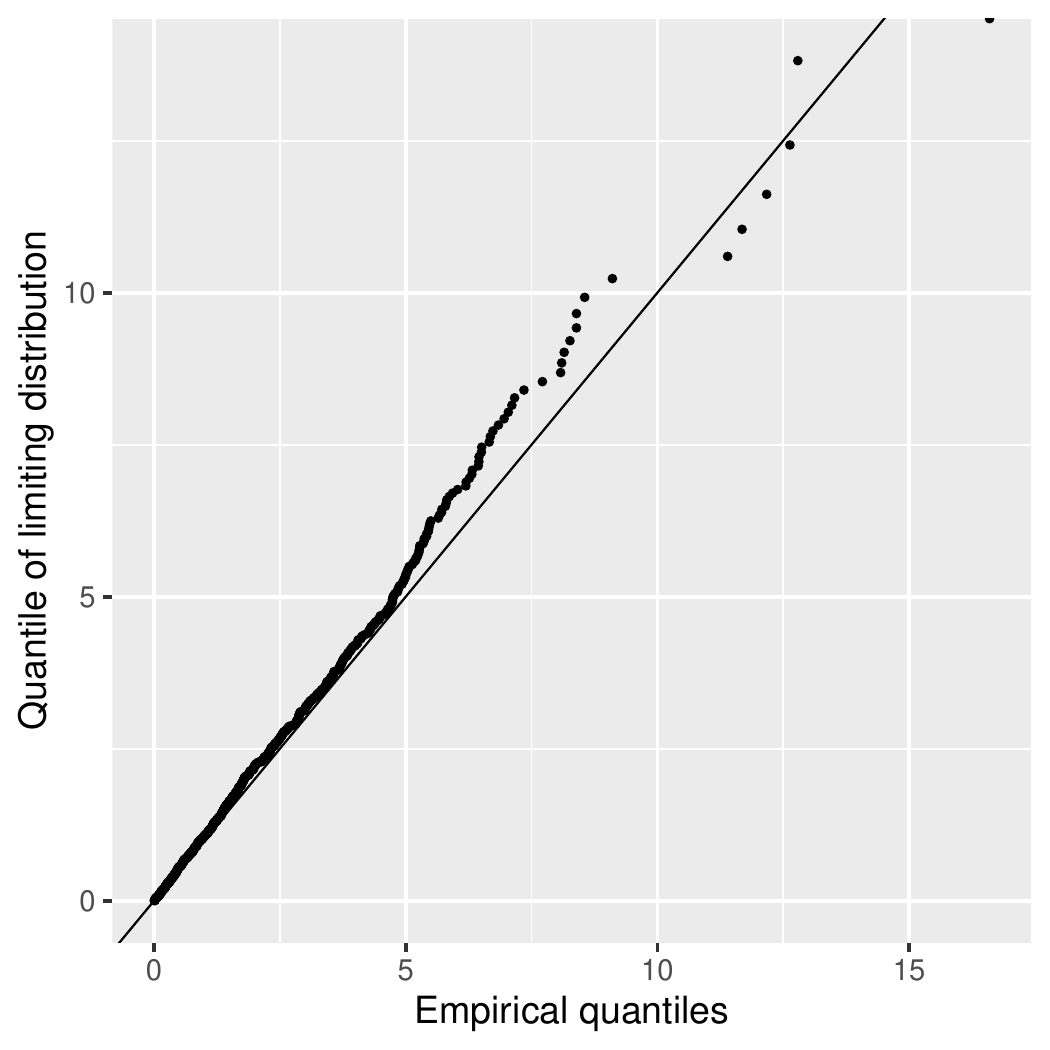}
		\caption*{}
	\end{subfigure}
	\begin{subfigure}[t]{0.48\textwidth}
		\includegraphics[width=6cm,height=6cm]{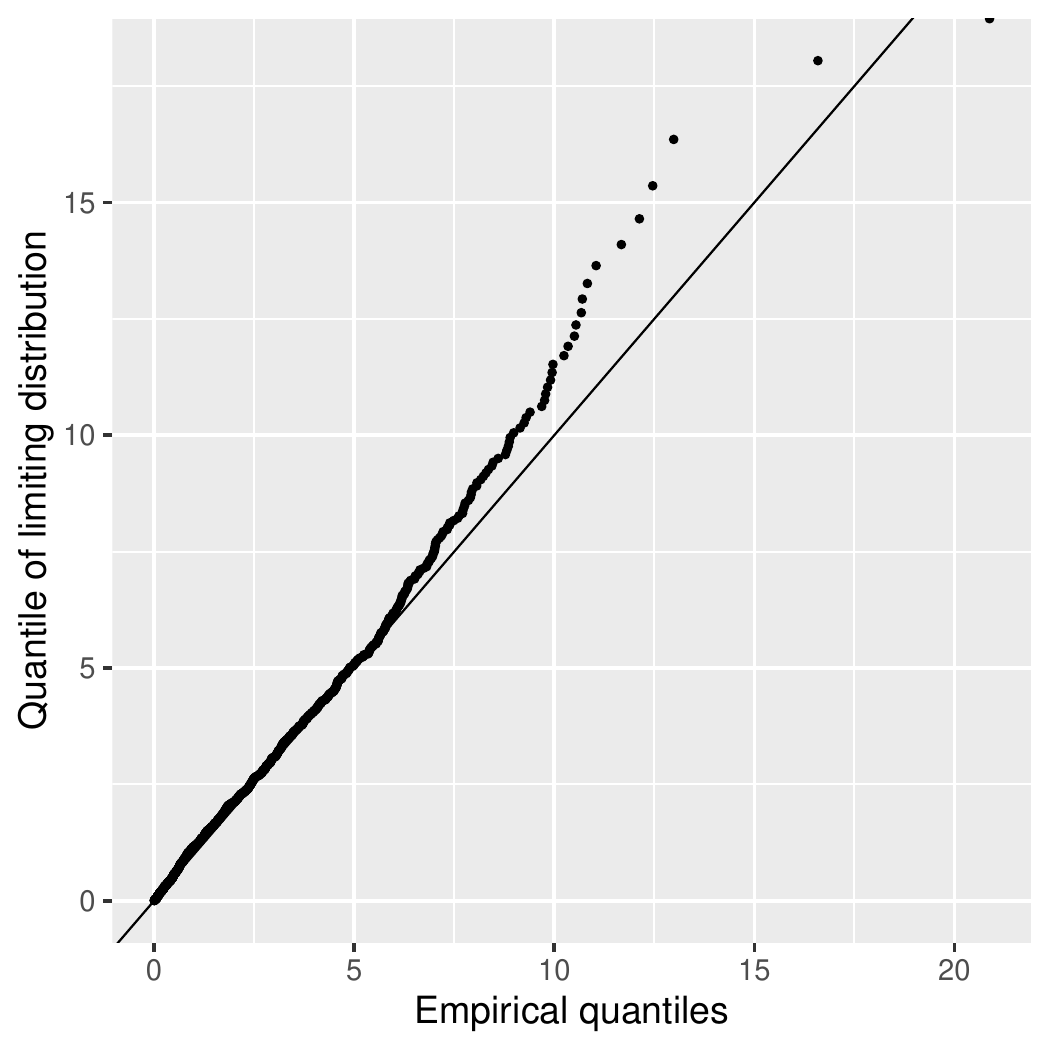}
		\caption*{}
	\end{subfigure}
	\caption{Left plot is Q-Q plot of 1000 replicates of $\tilde{T}_{\phi,n}(\theta_0)$,  $n=100$, under $H_0$; and $\phi(x)=-\log(x)+x-1$. Empirical distribution quantiles on the horizontal axis and the limiting $\chi^2_2$ distribution quantiles on the vertical axis. The right plot is corresponding plot under the alternative $\eta_n=\theta_0+\Delta n^{-1/2}$, with $\Delta=(1,1)$, and the limiting distribution $\chi^2_2\left(\frac{b_\phi^2}{\sigma_q^2}\Delta^t I(\theta_0)\Delta\right)$.}\label{fig4}
\end{figure}

\textit{Example 3.} Suppose that $X_1,\ldots,X_n$ is a random sample from a bivariate normal population $N(\eta,\Sigma)$, where 
$\eta= \begin{pmatrix}
	\mu_1 \\
	\mu_2 \\
\end{pmatrix}\in\mathbb{R}^2$ and $\Sigma= \begin{pmatrix}
	1 & 0 \\
	0 & 1 \\
\end{pmatrix}$. Consider testing the hypothesis $H_0:\eta=(0,0)^t$ against $H_A:\eta\neq(0,0)^t$. To examine the finite sample performance of the proposed test statistic $\tilde{T}_{\phi,n}(\theta_0)$, we take $n=100$ and $\phi(x)=-\log(x)+x-1$. Under $H_0$, the distribution of $\tilde{T}_{\phi,n}(\theta_0)$ is reasonably close the limiting $\chi^2_2$ distribution (see Figure \ref{fig4}). To assess performance of the test under local alternatives, we take $\Delta=(1,1)$ so that $\eta_n=\theta_0+\Delta n^{-1/2}=(0.1,0.1)^t$. Figure \ref{fig4} shows that under these local alternatives, the distribution of $\tilde{T}_{\phi,n}(\theta_0)$ is reasonably close to the limiting distribution $\chi^2_2\left(\dfrac{b_\phi^2}{\sigma_q^2}\Delta^t I(\theta_0)\Delta\right)$. 
\begin{figure}[!htb]
	\centering
	\begin{subfigure}[t]{0.48\textwidth}
		\includegraphics[width=6cm,height=6cm]{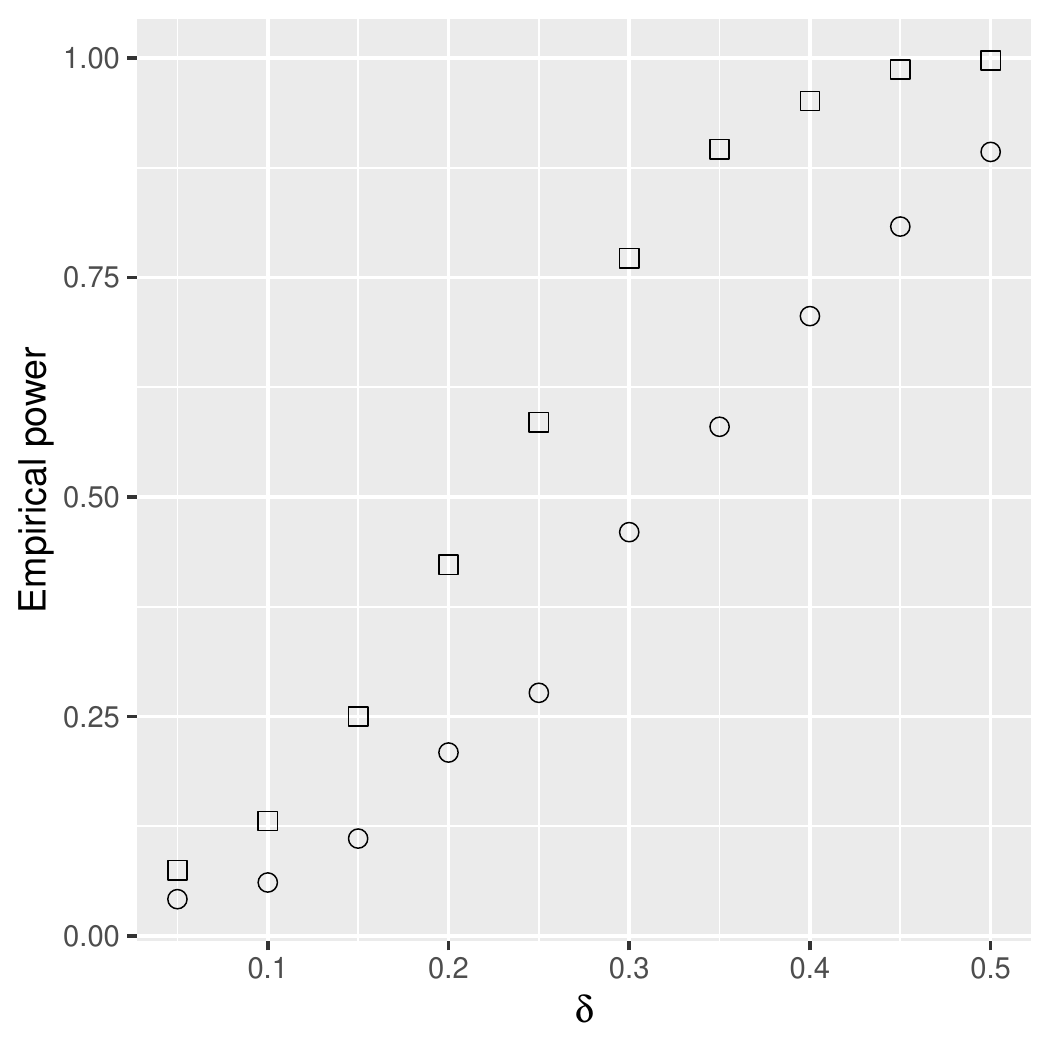}
		\caption*{{Alternative with $\eta=(0,\delta)$.} }
	\end{subfigure}
	\begin{subfigure}[t]{0.48\textwidth}
		\includegraphics[width=6cm,height=6cm]{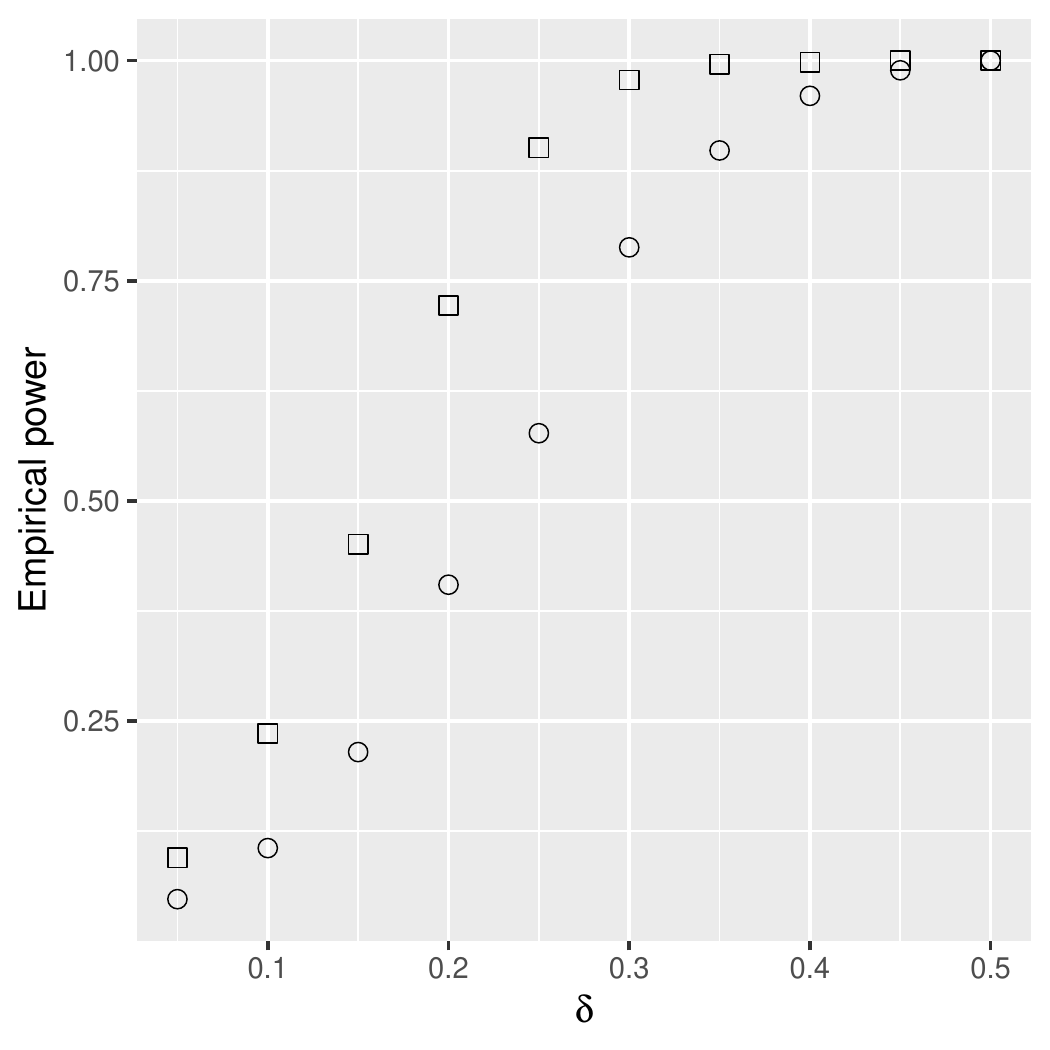}
		\caption*{{Alternative with $\eta=(\delta,\delta)$.}}
	\end{subfigure}
	\caption{{The empirical powers of $\tilde{T}_{\phi,n}(\theta_0)$, with $\phi(x)=-\log(x)+x-1$ (hollow circles) and the LRT (hollow squares), for $n=100$.}}\label{fig5}
\end{figure}
In this case, the LRT performs better than the proposed test (see Figure \ref{fig5}). For large differences, the performances are comparable. A reason for this could be that some features of normality are not captured by spacings.  {Performances of test statistics corresponding to some other choices of $\phi$ (e.g., $(x-1)\log(x)$, $(1-\sqrt{x})^2$) are similar to that of test statistic corresponding to $\phi(x)=-\log(x)+x-1$. The values of $b_\phi^2$ and $\sigma_q^2$ corresponding to these functions $\phi$ are given in Kuljus and Ranneby (2020).}  As the sample size increases, the agreement between observed quantiles of the test statistics and quantiles of the limiting distribution improves.
\FloatBarrier
\section{Real Data Analysis}
In this section, we use a real dataset to illustrate usefulness of the proposed methodology. Alzaatreh et al. (2012) analysed a fatigue life data of 6061-T6 aluminum coupons and found that a Gamma-Pareto distribution fits the data well. The distribution function of the Gamma-Pareto distribution (with parameters $\alpha,c,\beta>0$), is given by
\begin{align*}
	G_{\alpha,c,\beta}(x)=\begin{cases}
		\frac{1}{\Gamma(\alpha)}\Gamma(\alpha,c^{-1}\log(x/\beta)),~&\text{if }x>\beta\\
		0,~&\text{otherwise},
	\end{cases}
\end{align*}
where $\Gamma(a,x)=\int_{0}^{x}t^{a-1}e^{-t}~dt$, $x>0$, $a>0$, is the incomplete Gamma function. We consider the same data for analysis. This data consists of $101$ observations and it has maximum of $6$ ties. Since this data has ties, the asymptotically optimal parametric test based on simple spacings cannot be used and a test based on higher order spacings may be a remedy. We take $m=7$. For testing $H_0:\eta:=(\alpha,c,\beta)=\theta_0$ against $H_A:\eta\neq \theta_0$, for some $\theta_0$, the performance of the proposed tests are given in Table \ref{table1}. 
\begin{table}[h]
	\centering
	\caption{{For the fatigue life data of 6061-T6 aluminum coupons data, assuming Gamma-Pareto distribution, p-values of tests based on asymptotic distribution of $\tilde{T}_{-\log,n}^{(7)}(\theta_0)$.}}\label{table1}
	{\begin{tabular}{lcccc}
			\hline
			& \multicolumn{4}{c}{Case}                                  \\ \cline{2-5} 
			& I           & II            & III          & IV           \\ \hline 
			$\theta_0$ & $(15.0209,0.04258,70)$ & $(15.0209,0.04258,73)$ & $(15.0209,0.045,70)$ & $(16,0.04258,70)$ \\
			p-value & $0.5480513$ & $0.001349312$ & $0.01077901$ & $0.003335154$ \\ \hline
	\end{tabular}}
\end{table}
For the given data, the maximum likelihood estimator (MLE) of the parameters are $(\hat\alpha,\hat{c},\hat\beta)=(15.0209,0.04258,70)$. In the first case, we have considered the MLE as the $\theta_0$ and in other cases some deviations from the MLE have been considered. It appears from Table \ref{table1} that the proposed tests may be useful in some situations. 
\section{Discussion}
In this paper, we have studied several new parametric tests based on sample spacings, for testing simple and composite hypotheses concerning univariate and multivariate absolutely continuous distributions. 

For univariate distribution functions, under quite general conditions, the test based on the statistic $\tilde{T}_{\phi,n}^{(m)}(\theta_0)$  has similar asymptotic properties as the likelihood ratio test (LRT). For finite $m$, the test corresponding to $\phi(x)=-\log(x)$ gives an asymptotically efficient test. For $m\to\infty$, such that $m=o(n)$, there are multiple asymptotically efficient tests. For example, any test corresponding to $\phi=\phi_\gamma$ (see (\ref{pow_fun})) is asymptotically efficient. 

For multivariate distribution functions, spacings are defined in terms of the nearest neighbor balls. For simple multivariate spacings, the distribution of $\tilde{T}_{\phi,n}$ is similar to that of the LRT statistic. In this case we do not have any asymptotically efficient test. 

For moderately large sample size ($n\geq50$), our simulation study suggests that distributions of all the test statistics (corresponding to different convex functions $\phi$) are reasonably close to their limiting distribution under the null hypothesis as well as under local alternatives. An analysis based on real data also shows that small departures from the true distribution may be detected by the proposed tests. For certain distributions, the LRT does not exist (due to unboundedness of the likelihood function) and, in such situations, the proposed tests may be useful.

\section*{Appendix}
\begin{proof}[{Theorem \ref{thm2}}]
	For simplicity let us denote $S_{\phi,n}^{(m)}(\theta_0)$, $T_{\phi,n}^{(m)}(\theta_0)$ and $\hat{\theta}_{\phi,n}^{(m)}$ by $S$, $T$ and $\hat{\theta}$, respectively.
	Also, denote 
	\begin{align*}
		S_j'(\tau)=\frac{\partial}{\partial\theta_j}S(\theta)\bigg|_{\theta=\tau},\ 
		S_{j,k}''(\tau)=\frac{\partial^2}{\partial\theta_j\ \partial\theta_k}S(\theta)\bigg|_{\theta=\tau},\ 
		S_{j,k,l}'''(\tau)=\frac{\partial^3}{\partial\theta_j\ \partial\theta_k\ d\theta_l}S(\theta)\bigg|_{\theta=\tau},
	\end{align*}
	and $\delta=\theta_0-\hat{\theta}$, with $\delta_j$ as $j^{th}$ component of $\delta$.
	Using Taylor's expansion about $\hat{\theta}$, we have
	\begin{align}\label{app_eq1}
		T= \sum_{j=1}^{p}\delta_j S_j'(\hat{\theta})
		+\frac{1}{2}\sum_{j=1}^{p}\sum_{k=1}^{p}\delta_j\delta_k \left(S_{jk}''(\hat{\theta})+\frac{1}{3}\sum_{l=1}^{p}\delta_l S_{jkl}'''({\theta}^*)\right),
	\end{align}
	where $\theta^*$ is a point that lies on the line segment joining $\theta_0$ and $\hat{\theta}$. Under $H_0$, we have $\hat{\theta}\stackrel{p}{\rightarrow}\theta_{0}$. This implies the probability that $\hat{\theta}$ lies in any open neighbourhood of $\theta_0$ in $\Theta$ converges to 1, as $n\to\infty$. By the hypothesis of the theorem, we have $S_j'(\hat{\theta})=0;\ j=1,\ldots,p$. Observe that, using Taylor's expansion about $\theta_0$, we have
	\begin{align}\label{app_eq2}
		S_{jk}''(\hat{\theta})= S_{jk}''({\theta}_0)-\sum_{l=1}^{p}\delta_lS_{jkl}'''({\theta}^{**}),
	\end{align}
	where $\theta^{**}$ is a point on the line segment joining $\theta_0$ and $\hat{\theta}$. Using arguments similar to those used in the proof of Theorem 2, Ekstr\"om et al. (2020), we have, for $j,k=1,2,\ldots,p$, 
	\begin{align}\label{app_eq3}
		S_{jk}''({\theta}_0)+\frac{1}{3}\sum_{l=1}^{p}\delta_lS_{jkl}'''(\theta^*)-\sum_{l=1}^{p}\delta_lS_{jkl}'''({\theta}^{**})= S_{jk}''({\theta}_0)+o_p(1) \stackrel{p}{\rightarrow}E\left({\zeta}_{m}^{2} \phi ''({\zeta}_{m} )\right)I_{j,k}(\theta_0).
	\end{align}
	From (\ref{app_eq1})-(\ref{app_eq3}) and Theorem \ref{thm1}, we get
	\begin{align*}
		2nT&=E\left({\zeta}_{m}^{2} \phi ''({\zeta}_{m} )\right)n\sum_{j=1}^{p}\sum_{k=1}^{p}\delta_j\delta_kI_{jk}(\theta_0)+o_p(1)\\
		&=E\left({\zeta}_{m}^{2} \phi ''({\zeta}_{m} )\right) n(\theta_0-\hat{\theta})^t I(\theta_0)(\theta_0-\hat{\theta}) +o_p(1).
	\end{align*}
	Consequently, under $H_0$,
	\begin{align*}
		\tilde{T}_{\phi,n}^{(m)}=\frac{2nT}{E\left({\zeta}_{m}^{2} \phi ''({\zeta}_{m} )\right)\sigma_{\phi,m}^2}
		=\frac{n}{\sigma_{\phi,m}^2}(\theta_0-\hat{\theta})^t I(\theta_0)(\theta_0-\hat{\theta}) +o_p(1)\stackrel{d}{\rightarrow}\chi_p^2\text{ as }n\to\infty,
	\end{align*} 
	establishing part $(a)$ of the theorem.
	
	For $\eta_n=\theta_0+\Delta n^{-1/2}$, observe that
	\begin{align*}
		\sqrt{n}(\hat{\theta}-\theta_0)=\sqrt{n}(\hat{\theta}-{\eta_n})+\Delta
		\stackrel{d}{\rightarrow}N\left(\Delta,\ {\sigma_{\phi,m}^2} {I(\theta_{0})}^{-1}\right)\text{ as }n\to\infty.
	\end{align*}
	Using arguments similar to those used in part $(a)$, we get
	\begin{align*}
		\tilde{T}_{\phi,n}\stackrel{d}{\rightarrow}\chi_p^2(\sigma_{\phi,m}^{-2}\Delta^tI(\theta_0)\Delta)\text{ as }n\to\infty.
	\end{align*} 
	This completes the proof of part $(b)$.
\end{proof}
Theorem \ref{thm3} can be proved using arguments similar to those used in proving Theorem \ref{thm2} and part (b) of Theorem \ref{thm1}.
\begin{proof}[{Corollary \ref{cor2}}]
	From Theorem \ref{thm1} and using arguments similar to those in proving Corollary 1 of Ekstr\"om et.al. (2020), we have
	\begin{align*}
		\sqrt{n}({\hat{\theta}}_{\phi_2,n}^{(m)}-\theta_{0}) \stackrel{d}{\rightarrow} 
		N\left(0,\  {I(\theta_{0})}^{-1}\right)\text{ as }n\to\infty.
	\end{align*}
	 {Here, $\frac{\partial}{\partial \theta_j}S_{\phi_1,n}^{(m)}(\hat{\theta}_{\phi_2,n}^{(m)})\neq0$, however,  $\frac{\partial}{\partial \theta_j}S_{\phi_1,n}^{(m)}(\hat{\theta}_{\phi_2,n}^{(m)})-\frac{\partial}{\partial \theta_j}S_{\phi_1,n}^{(m)}(\hat{\theta}_{\phi_1,n}^{(m)})\xrightarrow{p}0$ as $n\to\infty$.} Thus using arguments similar to those in the proof of Theorem \ref{thm2}, we have under~$H_0$
	\begin{align*}
		{\tilde{T}}^{(m)}_{\phi_1,\phi_2,n}=&\dfrac{2n\left({S_{\phi_1,n}^{(m)}(\theta_0) -S_{\phi_1,n}^{(m)}(\hat{\theta}_{\phi_2,n}^{(m)})}\right)} {E\left({\zeta}_{m}^{2} \phi_1 ''({\zeta}_{m} )\right)\sigma_{\phi_2,m}^2}+o_p(1)\\
		=& \frac{n}{\sigma_{\phi_2,m}^2}(\theta_0-{\hat{\theta}_{\phi_2,n}^{(m)}})^t I(\theta_0)(\theta_0-{\hat{\theta}_{\phi_2,n}^{(m)}}) +o_p(1)
		\stackrel{d}{\rightarrow}\chi_p^2\text{ as }n\to\infty.
	\end{align*}
	Similarly, for $\eta_n=\theta_{0}+\Delta n^{-1/2}$, 
	\begin{align*}
		\sqrt{n}({{\hat{\theta}}_{\phi_2,n}^{(m)}}-\theta_{0}) \stackrel{d}{\rightarrow} 
		N\left(\Delta,\  {I(\theta_{0})}^{-1}\right)\text{ as }n\to\infty.
	\end{align*}
	Consequently,
	\begin{align*}
		{\tilde{T}}^{(m)}_{\phi_1,\phi_2,n}= 
		{n}(\theta_0-{\hat{\theta}_{\phi_2,n}^{(m)}})^t I(\theta_0)(\theta_0-{\hat{\theta}_{\phi_2,n}^{(m)}}) +o_p(1)\stackrel{d}{\rightarrow}\chi_p^2(\Delta^tI(\theta_{0})\Delta)\text{ as }n\to\infty.
	\end{align*}
\end{proof}
Corollary \ref{cor1} can be proved using arguments similar to those used in proving Corollary \ref{cor2}.
\begin{proof}[{Theorem \ref{thm4}}] The testing problem can be equivalently written as 
	\begin{align}\label{comp_hyp_eq}
		H_0:\eta=g(\beta)\ \text{ against }\ H_A:\eta\neq g(\beta),
	\end{align} 
	where $g=(g_1,g_2,\ldots,g_p):\mathbb{R}^{p-r}\to\mathbb{R}^p$ is a vector valued function such that  the $p\times (p-r)$ matrix $G(\beta)=((\frac{\partial}{\partial \beta_j}g_i(\beta)))$ exists and is continuous in $\beta$, and $rank(G(\beta))=p-r$. For example, let $p=3$, $\eta=(\eta_1,\eta_2,\eta_3)^t$ and $h(\eta)=\eta_1-\eta_2$. Then $\beta=(\beta_1,\beta_2)^t$ and $g(\beta)=(\beta_1,\beta_1,\beta_2)^t$. For the testing problem (\ref{comp_hyp_eq}), the proof is similar to the proof of Theorem 5.6.3 of Sen and Singer (1994). Details of the proof are given in the supplementary material.
\end{proof}
\begin{proof}[{Theorem \ref{thm6}}]
	Using Kuljus and Ranneby (2020), the proof follows along the lines of the proof of Theorem \ref{thm2}. Details of the proof are given in the supplementary material.
\end{proof}
\begin{proof}[{Theorem \ref{thm7}}]
	Using Kuljus and Ranneby (2020), the result follows on following the lines of the proof of Theorem \ref{thm4}. 
\end{proof}

\newpage

\title{Supplementary material to ``Some parametric tests based on sample spacings''}

\maketitle
\section*{Proof of Theorem 1}
The proof of Theorem  1, is based on Lemma 5.2, Lehman \& Casella (1998) and Theorem 2, Ekstr\"om et al. (2020).  First we state Lemma 5.2, Lehman \& Casella (1998), which is as follows.
\begin{lemma}
	Let $(T_{1n},T_{2n},\ldots,T_{pn})$ converges weakly to 
	$(T_{1},T_{2},\ldots,T_{p})$. Suppose for each $j$ and $k$, $A_{jkn}\xrightarrow{p}a_{jk}$, and $A=((a_{jk}))$ is non-singular. Let $B=((b_{jk}))=A^{-1}$. Then, if the distribution of $(T_{1},T_{2},\ldots,T_{p})$ has a density with respect to Lebesgue measure on $\mathbb{R}^p$, the solution of $(Y_{1n},Y_{2n},\ldots,Y_{pn})$ of
	\begin{align} \label{supp1:eq1}
		\sum_{k=1}^p A_{jkn}Y_{kn}=T_{jn},~j=1,2\ldots,p,
	\end{align}
	tend in probability to the solutions $(Y_{1},Y_{2},\ldots,Y_{p})$ of
	\begin{align} \label{supp1:eq2}
		\sum_{k=1}^p a_{jk}Y_{k}=T_{j},~j=1,2\ldots,p,
	\end{align}
	given by
	\begin{align} \label{supp1:eq3}
		Y_j=\sum_{k=1}^pb_{jk}T_k.
	\end{align}
\end{lemma}

\begin{proof}[\textcolor{blue}{Proof of Theorem 1}]
	Let $\theta_{0}$ be the true parameter. For simplicity let us denote $S_{\phi,n}^{(m)}(\theta)$ and $\hat{\theta}_{\phi,n}^{(m)}$, respectively, by $S$ and $\hat{\theta}$. Let $\theta=(\theta_{1},\theta_{2},\ldots,\theta_{p})^T$, $\theta_0=(\theta_{01},\theta_{02},\ldots,\theta_{0p})^T$ and $\hat\theta=(\hat\theta_{1},\hat\theta_{2},\ldots,\hat\theta_{p})^T$. 
	Also denote 
	\begin{align*}
		S_j'(\tau)=\frac{\partial}{\partial\theta_j}S(\theta)\bigg|_{\theta=\tau},\ 
		S_{j,k}''(\tau)=\frac{\partial^2}{\partial\theta_j\ \partial\theta_k}S(\theta)\bigg|_{\theta=\tau},\ 
		S_{j,k,l}'''(\tau)=\frac{\partial^3}{\partial\theta_j\ \partial\theta_k\ d\theta_l}S(\theta)\bigg|_{\theta=\tau}.
	\end{align*}
	
	\subsubsection*{(Consistency)} Let $B_\epsilon(\theta_{0})= \{\theta\in\Theta:||\theta-\theta_{0}||<\epsilon\}$  for $\epsilon>0$. 
	To prove the existence, with probability approaching 1, of a sequence of the GOSE, we shall examine behaviour of $S(\theta)$ on the sphere $B_\epsilon(\theta_{0})$. It is required to show that for any sufficiently small $\epsilon$,
	\begin{align} \label{supp1:eq4}
		P(S(\theta)>S(\theta_{0}))\xrightarrow{n\to\infty}1,~\forall \theta\in B_\epsilon(\theta_{0}).
	\end{align}
	Further, the proof is analogous to the one-dimensional case. The proof of ( \ref{supp1:eq4}) is similar to that of part (a) of Theorem 5.1 in Lehman and Casella (1998, p.~463), with some modifications.
	
	Observe that for sufficiently small $\epsilon$ and $\theta\in B_\epsilon(\theta_{0})$,
	\begin{align} \label{supp1:eq5}
		&S(\theta)-S(\theta_{0})\nonumber\\
		\geq&\sum_{j=1}^{p}(\theta_j-\theta_{0j})S_j'(\theta_0)+ \sum_{j=1}^{p} \sum_{k=1}^{p} (\theta_j-\theta_{0j})(\theta_k-\theta_{0k})S_{jk}''(\theta_0)\nonumber\\
		&- \sum_{j=1}^{p} \sum_{k=1}^{p} \sum_{l=1}^{p} (\theta_j-\theta_{0j})(\theta_k-\theta_{0k}) (\theta_l-\theta_{0l}) \frac{1}{M}\sum_{i=1}^MM_{jkl}\left(\frac{n+1}{m}D_{i:n}^{(m)}(\theta_{0})\right).
	\end{align}
	For each $j=1,2,\ldots,p$ and large $n$, with large probability $S_j'(\theta_0)$ is close to $S_j'(\hat\theta_{0j})$, where $\hat\theta_{0j}= (\theta_{01},\theta_{02},\ldots,\hat{\theta}_j,\ldots,\theta_{0p})^T$ (cf. Ekstr\"om et al. 2020). So, $S_j'(\theta_0)$ converges in probability to zero and hence the first term in the RHS of ( \ref{supp1:eq5}) converges in probability to zero. Using arguments similar to those in  Ekstr\"om et al. (2020), we can show that $S_{jk}''(\theta_0)$ converges in probability to $\mathbb{E}(\zeta_m^2\phi''(\zeta_m))I_{jk}(\theta_0)$. So, the second term in the RHS of ( \ref{supp1:eq5}) converges in probability to a quadratic form. Similarly, we can show that the third term in the RHS of ( \ref{supp1:eq5}) converges in probability to zero. Thus, for sufficiently small $\epsilon$ and $\theta\in B_\epsilon(\theta_{0})$, $S(\theta)-S(\theta_{0})$ is non-negative with probability tending to one. This proves the statement ( \ref{supp1:eq4}).
	\subsubsection*{(Asymptotic normality)}
	Under assumptions of Theorem 1, we have
	\begin{align*}
		S'_j(\hat{\theta})=0,~\forall j=1,2,\ldots,p.
	\end{align*}
	Using the Taylor expansion, we get
	\begin{align*}
		0=S_j'(\hat{\theta})=S'_j (\theta _{0})+\sum_{k=1}^{p}(\hat{\theta }_k -\theta _{0k} )S''_{jk}(\theta _{0})+\frac{1}{2}\sum_{k=1}^{p}\sum_{l=1}^{p} (\hat{\theta }_{k} -\theta _{0k} )(\hat{\theta }_{l} -\theta _{0l} ) S'''_{jkl} (\tilde{\theta }_{j}),
	\end{align*}
	where $\tilde{\theta }_{j}$ lies on the line joining $\hat{\theta }$ and $\theta _{0}$.
	
	It follows that
	\begin{align*}
		\sqrt{n}\sum_{k=1}^{p}(\hat{\theta }_k -\theta _{0k} )\left[S''_{jk}(\theta _{0})+\frac{1}{2}\sum_{l=1}^{p}(\hat{\theta }_{l} -\theta _{0l} ) S'''_{jkl} (\tilde{\theta }_{j})\right]=-\sqrt{n}S'_j (\theta _{0}).
	\end{align*}
	This has the form of  \ref{supp1:eq1} with
	\begin{align*}
		Y_{kn}=&\sqrt{n}(\hat{\theta }_k -\theta _{0k} ),\\
		A_{jkn}=&S''_{jk}(\theta _{0})+\frac{1}{2}\sum_{l=1}^{p}(\hat{\theta }_{l} -\theta _{0l} ) S'''_{jkl} (\tilde{\theta }_{j}),\\
		T_{jn}=& -\sqrt{n}S'_j (\theta _{0}).
	\end{align*}
	Using similar arguments as in the proof of Theorem 2, Ekstr\"om et al. (2020), every linear combination of $(T_{1n},T_{2n},\ldots,T_{pn})$ has limiting normal distribution. Thus, using Cramer-Wald device
	\begin{align*}
		(T_{1n},T_{2n},\ldots,T_{pn})\xrightarrow{d}N(0,\sigma_{\phi,m}^2(\mathbb{E}(\zeta_m^2\phi''(\zeta_m)))^2I(\theta_0)).
	\end{align*}
	Further using arguments as in the proof of Theorem 2, Ekstr\"om et al. (2020),
	\begin{align*}
		A_{jkn}\xrightarrow{p}a_{jk}=~\mathbb{E}(\zeta_m^2\phi''(\zeta_m))I_{jk}(\theta_0).
	\end{align*}
	Thus, the limiting distribution of $(Y_{1n},Y_{2n},\ldots,Y_{pn})$ is the distribution of $(Y_{1},Y_{2},\ldots,Y_{p})$. The $(Y_{1},Y_{2},\ldots,Y_{p})$ are given by solution of
	\begin{align*}
		&\sum_{k=1}^p \mathbb{E}(\zeta_m^2\phi''(\zeta_m))I_{jk}(\theta_{0})Y_{k}=T_{j},~j=1,2\ldots,p,\\
		\Rightarrow&~ (Y_1,Y_2,\ldots,Y_p)^t=~ \frac{1}{\mathbb{E}(\zeta_m^2\phi''(\zeta_m))} [I(\theta_{0})]^{-1}(T_1,T_2,\ldots,T_p)^t.
	\end{align*}
	This concludes the proof.
\end{proof}
\section*{Proof of Theorem 4}
The proof of Theorem  4, is similar to that of Theorem 5.6.3 of Sen \& Singer (1994) with some modifications. 
\begin{proof}[\textcolor{blue}{Proof of Theorem 4}] The testing problem can be equivalently written as 
	\begin{align*}
		H_0:\eta=g(\beta)\ \text{ against }\ H_A:\eta\neq g(\beta),
	\end{align*} 
	where $g=(g_1,g_2,\ldots,g_p):\mathbb{R}^{p-r}\to\mathbb{R}^p$ is a vector valued function such that  the $p\times (p-r)$ matrix $G(\beta)=((\frac{\partial}{\partial \beta_j}g_i(\beta)))$ exists and is continuous in $\beta$, and $rank(G(\beta))=p-r$. For example, let $p=3$, $\eta=(\eta_1,\eta_2,\eta_3)^t$, $h(\eta)=\eta_1-\eta_2$ then $\beta=(\beta_1,\beta_2)^t$ and $g(\beta)=(\beta_1,\beta_1,\beta_2)^t$. 
	
	Note that if $\hat{\beta}_n$ is GSE of $\beta$, it follows
	\begin{align*}
		\inf_{\beta\in\mathbb{R}^{p-r}}S_{\phi,n}^{(m)}(\beta)=
		\inf_{\{\theta\in\Theta:h(\theta)=0\}}S_{\phi,n}^{(m)}(\theta).
	\end{align*}
	Thus the test statistic is 
	\begin{align}
		T_{\phi,n}^{(m)}(h)=&\frac{2n}{E\left({\zeta}_{m}^{2} \phi ''({\zeta}_{m} )\right)\sigma_{\phi,m}^2}[S_{\phi,n}^{(m)}(\hat\beta_n)- S_{\phi,n}^{(m)}(\hat\theta_n)]\nonumber\\
		=&\frac{2n}{\mathbb{E}\left({\zeta}_{m}^{2} \phi ''({\zeta}_{m} )\right)\sigma_{\phi,m}^2} [S_{\phi,n}^{(m)}(\theta_0)-S_{\phi,n}^{(m)}(\hat\theta_n)]\nonumber\\
		&- \frac{2n}{\mathbb{E}\left({\zeta}_{m}^{2} \phi ''({\zeta}_{m} )\right)\sigma_{\phi,m}^2} [S_{\phi,n}^{(m)}(\theta_0)-S_{\phi,n}^{(m)}(\hat\beta_n)]\nonumber\\
		=& Q_1-Q_2, \text{ say.}
	\end{align}
	Using Taylor's expansion, we have 
	\begin{align}
		0=&\sqrt{n}\frac{\partial}{\partial \theta} S_{\phi,n}^{(m)}(\theta)\bigg|_{\theta=\hat{\theta}_n}\nonumber\\
		=& \sqrt{n} S_{\phi,n}^{(m)'}(\theta_0)+ \frac{\partial^2}{\partial \theta\partial \theta^t} S_{\phi,n}^{(m)}(\theta)\bigg|_{\theta=\tilde{\theta}_n}\sqrt{n}(\hat{\theta}_n-\theta_0)
		,
	\end{align}
	where $\tilde{\theta}_n$ is on the line segment joining $\theta_0$ and $\hat{\theta}_n$. Using arguments similar to those in the proof of Theorem 2, Ekstr\"om et al. (2020), we obtain
	\begin{align}
		&\frac{\partial^2}{\partial \theta\partial \theta^t} S_{\phi,n}^{(m)}(\theta)\bigg|_{\theta=\tilde{\theta}_n} \xrightarrow{p}~ \mathbb{E}\left({\zeta}_{m}^{2} \phi ''({\zeta}_{m} )\right)I(\theta_0),\text{ as 
		}n\to\infty,\nonumber\\
		\Rightarrow& \sqrt{n}(\hat{\theta}_n-\theta_0)=\frac{-I(\theta_0)^{-1}}{\mathbb{E}\left({\zeta}_{m}^{2} \phi ''({\zeta}_{m} )\right)} \sqrt{n}S_{\phi,n}^{(m)'}(\theta_0)+o_p(1).
	\end{align}
	Similarly, we obtain
	\begin{align}
		\sqrt{n}(\hat{\beta}_n-\beta_0)=\frac{-I^*(\beta_0)^{-1}}{\mathbb{E}\left({\zeta}_{m}^{2} \phi ''({\zeta}_{m} )\right)} \sqrt{n}S_{\phi,n}^{(m)'}(\beta_0)+o_p(1).
	\end{align}
	Using arguments similar to those in Theorem 2 and above relation, we get 
	\begin{align*}
		Q_1=&~\frac{2n}{\mathbb{E}\left({\zeta}_{m}^{2} \phi ''({\zeta}_{m} )\right)\sigma_{\phi,m}^2}~ [S_{\phi,n}^{(m)}(\theta_0)-S_{\phi,n}^{(m)}(\hat\theta_n)]\\
		=&~ \frac{n}{\sigma_{\phi,m}^2}~(\hat{\theta}_n-\theta_0)^t ~I(\theta_0)~(\hat{\theta}_n-\theta_0)+o_p(1)\\
		=&~ \frac{n}{\sigma_{\phi,m}^2[\mathbb{E}\left({\zeta}_{m}^{2} \phi ''({\zeta}_{m} )\right)]^2}~[S_{\phi,n}^{(m)'}(\theta_0)]^t~I(\theta_0)^{-1}~S_{\phi,n}^{(m)'}(\theta_0) +o_p(1).
	\end{align*}
	Similarly, we obtain
	\begin{align*}
		Q_2	=&~ \frac{n}{\sigma_{\phi,m}^2[\mathbb{E}\left({\zeta}_{m}^{2} \phi ''({\zeta}_{m} )\right)]^2}~[S_{\phi,n}^{(m)'}(\beta_0)]^t~I^*(\beta_0)^{-1}~S_{\phi,n}^{(m)'}(\beta_0) +o_p(1).
	\end{align*}
	Observe that
	\begin{align*}
		S_{\phi,n}^{(m)'}(\beta)=&\frac{\partial}{\partial\beta}~S_{\phi,n}^{(m)}(\beta)= 
		\frac{\partial}{\partial\beta}~g(\beta)\frac{\partial}{\partial g(\beta)}~S_{\phi,n}^{(m)}(g(\beta))\nonumber\\
		=&[G(\beta)]^tS_{\phi,n}^{(m)'}(\theta).
	\end{align*}
	Now we have 
	\begin{align*}
		&\sqrt{n}S_{\phi,n}^{(m)'}(\beta_0)\xrightarrow{d}~N(0,\sigma_{\phi,m}^2[\mathbb{E}\left({\zeta}_{m}^{2} \phi ''({\zeta}_{m} )\right)]^2I^*(\beta_0))\\
		\text{and }& \sqrt{n}[G(\beta)]^tS_{\phi,n}^{(m)'}(\theta_0)\xrightarrow{d}~N(0,\sigma_{\phi,m}^2[\mathbb{E}\left({\zeta}_{m}^{2} \phi ''({\zeta}_{m} )\right)]^2[G(\beta)]^tI(\theta_0)[G(\beta)]).
	\end{align*}
	Thus,
	\begin{align*}
		I^*(\beta_0)=[G(\beta)]^tI(\theta_0)[G(\beta)]
	\end{align*}
	and
	\begin{align*}
		Q_2	=&~ \frac{n}{\sigma_{\phi,m}^2[\mathbb{E}\left({\zeta}_{m}^{2} \phi ''({\zeta}_{m} )\right)]^2}~[S_{\phi,n}^{(m)'}(\theta_0)]^t~[G(\beta)]I^*(\beta_0)^{-1}[G(\beta)]^t~S_{\phi,n}^{(m)'}(\theta_0) +o_p(1).
	\end{align*}
	Hence,
	\begin{align*}
		&T_{\phi,n}^{(m)}(h)\nonumber\\
		=&~\frac{n}{\sigma_{\phi,m}^2[\mathbb{E}\left({\zeta}_{m}^{2} \phi ''({\zeta}_{m} )\right)]^2}~[S_{\phi,n}^{(m)'}(\theta_0)]^t~\big[I(\theta_0)^{-1}-[G(\beta)]I^*(\beta_0)^{-1}[G(\beta)]^t\big]~S_{\phi,n}^{(m)'}(\theta_0) +o_p(1).
	\end{align*}
	Now, we have
	\begin{align*}
		&\sqrt{n}S_{\phi,n}^{(m)'}(\theta_0)\xrightarrow{d}N(0,\sigma_{\phi,m}^2[\mathbb{E}\left({\zeta}_{m}^{2} \phi ''({\zeta}_{m} )\right)]^2I(\theta_0))\\
		\text{and }& \big[I(\theta_0)^{-1}-[G(\beta)]I^*(\beta_0)^{-1}[G(\beta)]^t\big]I(\theta_0)\text{ is idempotent.}
	\end{align*}
	Next,
	\begin{align*}
		tr\bigg[\big[I(\theta_0)^{-1}-[G(\beta)]I^*(\beta_0)^{-1}[G(\beta)]^t\big]I(\theta_0) \bigg]=r.
	\end{align*}
	Thus,
	\begin{align*}
		T_{\phi,n}^{(m)}(h)\xrightarrow{d}~\chi^2_{(r)}\text{ as }n\to\infty.
	\end{align*}
	This concludes the proof.
\end{proof}
\section*{Proof of Theorem 5}
Using Kuljus and Ranneby (2020), the proof follows on the lines of the proof of Theorem 2. Detail of the proof is given below.
\begin{proof}
	For simplicity, let us denote $S_{\phi,n}(\theta)$ and $\hat{\theta}_{\phi,n}$ by $S$ and $\hat{\theta}$, respectively. Let $\theta=(\theta_{1},\theta_{2},\ldots,\theta_{p})^T$, $\theta_0=(\theta_{01},\theta_{02},\ldots,\theta_{0p})^T$ and $\hat\theta=(\hat\theta_{1},\hat\theta_{2},\ldots,\hat\theta_{p})^T$. 
	Also denote 
	\begin{align*}
		S_j'(\tau)=\frac{\partial}{\partial\theta_j}S(\theta)\bigg|_{\theta=\tau},\ 
		S_{j,k}''(\tau)=\frac{\partial^2}{\partial\theta_j\ \partial\theta_k}S(\theta)\bigg|_{\theta=\tau},
	\end{align*}
	and $\delta=\theta_0-\hat{\theta}$, with $\delta_j$ as $j^{th}$ component of $\delta$.
	Using Taylor's expansion about $\hat{\theta}$, we have
	\begin{align*}
		T_{\phi,n}(\theta_0)= \sum_{j=1}^{p}\delta_j S_j'(\hat{\theta})
		+\frac{1}{2}\sum_{j=1}^{p}\sum_{k=1}^{p}\delta_j\delta_k S_{jk}''({\theta}^*),
	\end{align*}
	where $\theta^*$ is a point that lies on the line segment joining $\theta_0$ and $\hat{\theta}$. Under $H_0$, we have $\hat{\theta}\stackrel{p}{\rightarrow}\theta_{0}$. This implies the probability that $\hat{\theta}$ lies in any open neighbourhood of $\theta_0$ in $\Theta$ converges to 1, as $n\to\infty$. Using the fact that $S_j'(\hat{\theta})=0;\ j=1,\ldots,p$, we have
	\begin{align*}
		T_{\phi,n}(\theta_0)
		&= \frac{1}{2}\sum_{j=1}^{p}\sum_{k=1}^{p}\delta_j\delta_k S_{jk}''({\theta}^*).
	\end{align*}
	To complete the proof, we need to show that $S_{jk}''({\theta}^*)$ converges in probability to $\mathbb{E}(Z_1^2\phi''(Z_1))I_{jk}(\theta_{0})$, where $I_{jk}(\theta)$ denotes the element of $j^{th}$ row and $k^{th}$ column of the matrix $I(\theta)$. Using the approach of Kuljus and Ranneby (2020, p. 984), we obtain $S_{jk}''({\theta}^*)$ converges in probability to $\mathbb{E}(Z_1^2\phi''(Z_1))I_{jk}(\theta_{0})$.
\end{proof}


\end{document}